\documentclass{amsproc}

\newtheorem{theorem}{Theorem}[section]
\newtheorem{lemma}[theorem]{Lemma}

\theoremstyle{definition}

\newtheorem{example}[theorem]{Example}

\theoremstyle{remark}
\newtheorem{remark}[theorem]{Remark}

\numberwithin{equation}{section}

\newcommand{\abs}[1]{\lvert#1\rvert}

\newcommand{\blankbox}[2]{%
  \parbox{\columnwidth}{\centering

    \setlength{\fboxsep}{0pt}%
    \fbox{\raisebox{0pt}[#2]{\hspace{#1}}}%
  }%
}
\newcommand{\R}{{\mathbb R}}
\newcommand{\D}{{\mathbb D}}
\newcommand{\HH}{{\mathbb H}}
\newcommand{\PP}{{\mathcal P}}
\newcommand{\boundary}{\partial_\infty\D}
\newcommand{\EE}{{\mathcal E\mathcal P}}
\newcommand{\ii}{{\rm Ind_f}}

\begin{document}

\title{Continua as minimal sets of homeomorphisms of $S^2$}

\author{Shigenori Matsumoto}
\address{Department of Mathematics, College of
Science and Technology, Nihon University, 1-8-14 Kanda, Surugadai,
Chiyoda-ku, Tokyo, 101-8308 Japan
}
\curraddr{Department of Mathematics, College of
Science and Technology, Nihon University, 1-8-14 Kanda, Surugadai,
Chiyoda-ku, Tokyo, 101-8308 Japan}
\email{matsumo@math.cst.nihon-u.ac.jp}

\author{Hiromichi Nakayama}
\address{Department of Physics and Mathematics, College of
Science and Engineering, Aoyama Gakuin University, 5-10-1 Huchinobe, Sagamihara, Kanagawa, Japan, 229-8558, Japan}
\curraddr{Department of Physics and Mathematics, College of
Science and Engineering, Aoyama Gakuin University, 5-10-1 Huchinobe, Sagamihara, Kanagawa, Japan, 229-8558, Japan}
\email{nakayama@gem.aoyama.ac.jp}
\thanks{Partially supported by Grant-in-Aid for Scientific Research (No.\ 20540096 for the first author and No.\ 19540090 for the second), Japan Society for the Promotion of Science, Japan}



\maketitle

\begin{abstract}
Let $f$ be an orientation preserving homeomorphism of $S^2$ which
has a  continuum $X$ as a minimal set. 
Then there are exactly two connected components of $S^2\setminus X$ 
which are left invariant by $f$ and all the others
are wandering. The Carath\'eodory rotation number of an invariant
component is irrational.
\end{abstract}

\section{Introduction}
Let $f$ be an orientation preserving homeomorphism of $S^2$ which
has a  continuum $X$ as a minimal set.  By a {\em continuum} we
mean a compact connected subset which is not a single point.
We have a great variety of examples of such homeomorphisms.
The simplest one is 
an irrational rotation on $S^2$, with a round circle as a 
minimal set. Besides this
a pathological diffeomorphism of $S^2$ is constructed in \cite{Ha}
which has a pseudo-circle as a minimal set. See also \cite{He} for
a curious diffeomorphisms.
Also a homeomorphism of $S^2$ with a minimal set homeomorphic
to a variant of the Warsaw circle is constructed in \cite{W}.
A fast approximation by conjugacy method is discussed in \cite{FK}, which
may produce such diffeomorphisms with various topological natures.

In all these examples the minimal sets $X$ separate $S^2$ into
two domains. So it is natural to ask if this is the case
with any minimal continuum.
It is well known that for any $n\in\mathbb N$, there is a continuum $X$ in
$S^2$ which separates $S^2$ into $n$ open domains $U_1,\cdots, U_n$
such that the frontier of each $U_i$ coincides with $X$ (\cite{K}). 

A connected
component $U$ of $S^2\setminus X$ is called an {\em invariant domain}
if $fU=U$, a {\em periodic domain} if $f^nU=U$ for some
$n\geq1$, and a {\em wandering domain} otherwise.

\begin{theorem} \label{0}
Consider an orientation preserving homeomorphism of $S^2$
which admits a continuum as a minimal set.
Then there are exactly two invariant domains and all the other domains
are wandering. The Carath\'eodory rotation numbers of both invariant domains are identical and irrational.
\end{theorem}

The overall strategy to prove Theorem \ref{0} is to use the
Carath\'eodory prime end theory and to apply the Cartwright-Littlewood
theorem.
Sections 2 and 3 are expositions of the prime end theory and the
Cartwright-Littlewood theorem, which are included since 
they are short and self-contained, and some
special features remarked in these sections are needed in 
the development of Section 4, which is devoted
to the proof of Theorem \ref{0}. Both Sections 2 and 3 concern
simply connected domains of closed oriented surfaces of any genus, and
Section 4 solely orientation preserving homeomorphisms of
the sphere $S^2$.
In Section 5 we will construct a homeomorphism which actually
admits a wandering domain.

\section{Prime ends}

Denote by $\Sigma$ a closed oriented surface
equipped with a smooth Riemannian metric $g$
and the associated area form $dvol$. Let $U\subset \Sigma$ be a {\em hyperbolic domain}
i.\ e.\ an open simply connected
subset such that $\Sigma\setminus U$ is not a singleton. 
(A nonhyperbolic simply connected domain exists only on the 2-sphere.)
The purpose of this section is to show that a homeomorphism of $U$
which extends to a homeomorphism of the closure $\overline U$
does extend to a homeomorphism of the so called Carath\'eodory
compactification $\hat U$, a closed disc. Here we are only concerned with
a simply connected domain in $\Sigma$. 
But there are generalizations to more general domains, which are found in
\cite{E} and \cite{M}. 
As general references of prime end theory, see also Sect.\ 17, \cite{Mi} and Chapter IX, \cite{T}.
The proof of the
main lemma here (Lemma \ref{1}) is taken from \cite{E}.

Let $0\in U$
be a base point.
A real line properly embedded in $U$ and not passing through $0$
is called a {\em cross cut}. 
A cross cut $c$ separates $U$ into two hyperbolic domains, as can
be seen by considering the one point compactification of $U$
and applying the Jordan curve theorem.
One of them 
not containing $0$ is called the {\em content} of $c$ and denoted by
$U(c)$. 
A sequence of cross cuts $\{c_i\}_{i=1}^\infty$ is called a {\em chain} if
$c_{i+1}\subset U(c_i)$ for each $i$.
Two chains $\{c_i\}$ and $\{c'_i\}$ are called {\em equivalent} if
for any $i$, there is $j$ such that $c'_j \subset U(c_i)$ and 
$c_j\subset U(c'_i)$.
An equivalence class of chains is called an {\em end} of $U$.
(This is quite different from the notion of ends for general noncompact spaces developed by H. Freudenthal et al., and exposed e.\ g.\ in
\cite{E2}.)
A homeomorphism between two hyperbolic domains
induces in an obvious way a bijection between the sets of ends.
Given an end $\xi$, the relatively closed set $C(\xi)=\cap_iU(c_i)$
is independent of the choice of a chain $\{c_i\}$ from the end $\xi$,
and is called the {\em content} of $\xi$.

\begin{figure}[htbp]
\begin{center}
\unitlength 0.1in
\begin{picture}( 30.0000, 20.0000)(  5.0000,-25.0000)
%
\special{pn 8}%
\special{pa 500 500}%
\special{pa 3500 500}%
\special{fp}%
%
\special{pn 8}%
\special{pa 500 500}%
\special{pa 500 2500}%
\special{fp}%
%
\special{pn 8}%
\special{pa 3500 500}%
\special{pa 3500 2500}%
\special{fp}%
%
\special{pn 13}%
\special{pa 2000 800}%
\special{pa 2000 2500}%
\special{fp}%
%
\special{pn 13}%
\special{pa 500 1220}%
\special{pa 1750 1220}%
\special{fp}%
%
\special{pn 8}%
\special{pa 1200 1000}%
\special{pa 1200 1220}%
\special{dt 0.045}%
%
\special{pn 13}%
\special{pa 500 1660}%
\special{pa 1750 1660}%
\special{fp}%
%
\special{pn 13}%
\special{pa 500 2000}%
\special{pa 1750 2000}%
\special{fp}%
%
\special{pn 13}%
\special{pa 500 2240}%
\special{pa 1750 2240}%
\special{fp}%
%
\special{pn 13}%
\special{pa 500 2380}%
\special{pa 1750 2380}%
\special{fp}%
%
\special{pn 8}%
\special{pa 500 2460}%
\special{pa 1750 2460}%
\special{fp}%
%
\special{pn 8}%
\special{pa 500 2500}%
\special{pa 3500 2500}%
\special{fp}%
%
\special{pn 8}%
\special{pa 1500 1220}%
\special{pa 1500 1450}%
\special{dt 0.045}%
%
\special{pn 8}%
\special{pa 1200 1450}%
\special{pa 1200 1660}%
\special{dt 0.045}%
%
\special{pn 8}%
\special{pa 1200 1840}%
\special{pa 1200 2000}%
\special{dt 0.045}%
%
\special{pn 8}%
\special{pa 1500 1660}%
\special{pa 1500 1840}%
\special{dt 0.045}%
%
\special{pn 8}%
\special{pa 1500 2000}%
\special{pa 1500 2150}%
\special{dt 0.045}%
%
\special{pn 8}%
\special{pa 1500 2240}%
\special{pa 1500 2330}%
\special{dt 0.045}%
%
\special{pn 8}%
\special{pa 1500 2380}%
\special{pa 1500 2420}%
\special{dt 0.045}%
%
\special{pn 8}%
\special{pa 1200 2150}%
\special{pa 1200 2240}%
\special{dt 0.045}%
%
\special{pn 8}%
\special{pa 1200 2330}%
\special{pa 1200 2380}%
\special{dt 0.045}%
%
\special{pn 8}%
\special{pa 1200 2420}%
\special{pa 1200 2460}%
\special{dt 0.045}%
%
\special{pn 13}%
\special{pa 1000 1000}%
\special{pa 2500 1000}%
\special{fp}%
%
\special{pn 13}%
\special{pa 1000 1450}%
\special{pa 2500 1450}%
\special{fp}%
%
\special{pn 13}%
\special{pa 1000 1840}%
\special{pa 2500 1840}%
\special{fp}%
%
\special{pn 13}%
\special{pa 1000 2150}%
\special{pa 2500 2150}%
\special{fp}%
%
\special{pn 13}%
\special{pa 1000 2330}%
\special{pa 2500 2330}%
\special{fp}%
%
\special{pn 13}%
\special{pa 1000 2420}%
\special{pa 2500 2420}%
\special{fp}%
%
\special{pn 8}%
\special{pa 2500 1000}%
\special{pa 3300 2500}%
\special{dt 0.045}%
%
\special{pn 8}%
\special{pa 2500 1450}%
\special{pa 3080 2500}%
\special{dt 0.045}%
%
\special{pn 8}%
\special{pa 2500 1840}%
\special{pa 2880 2500}%
\special{dt 0.045}%
%
\special{pn 8}%
\special{pa 2500 2150}%
\special{pa 2720 2500}%
\special{dt 0.045}%
%
\special{pn 8}%
\special{pa 2500 2330}%
\special{pa 2610 2500}%
\special{dt 0.045}%
%
\special{pn 8}%
\special{pa 2500 2420}%
\special{pa 2550 2500}%
\special{dt 0.045}%
\put(12.7000,-6.7000){\makebox(0,0)[lt]{$U$}}%
\put(11.7000,-11.6000){\makebox(0,0)[rb]{$c_1$}}%
\put(11.5000,-16.1000){\makebox(0,0)[rb]{$c_3$}}%
\put(15.4000,-13.9000){\makebox(0,0)[lb]{$c_2$}}%
\put(15.3000,-18.0000){\makebox(0,0)[lb]{$c_4$}}%
\put(28.2000,-15.4000){\makebox(0,0)[lb]{$c_1'$}}%
\put(27.4000,-18.8000){\makebox(0,0)[lb]{$c_2'$}}%
\put(26.9000,-22.0000){\makebox(0,0)[lb]{$c_3'$}}%
\end{picture}%
\caption{Topological chains}
\label{fig1}
\end{center}
\end{figure}
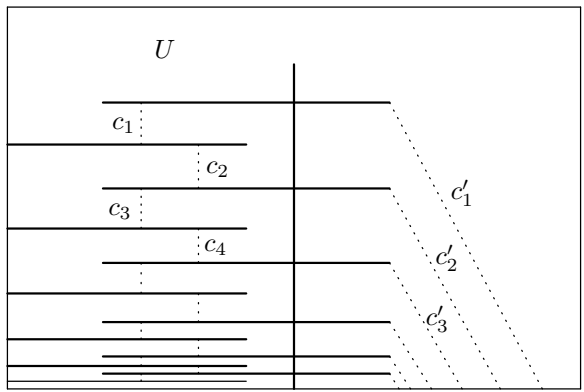

A chain
$\{c_i\}$  is called {\em topological} if
 the closures $\overline c_i$ of $c_i$ in $\Sigma$
are mutually disjoint and the diameter
${\rm diam}(c_i)$ converges to $0$ as $i\to\infty$.
Examples of topological chains, $\{c_i\}$ and $\{c'_i\}$,
 are given in Figure \ref{fig1}.
An end is called {\em prime} if it admits a topological chain.

\begin{lemma} \label{0a}
The content $C(\xi)$ of a prime end $\xi$ is empty.
\end{lemma}

{\bf Proof}. 
Assume the contrary and choose a point $x$ from $C(\xi)$.
Consider an arc $\gamma$ in $U$ joining $0$ to $x$. See Figure \ref{fig2}.
\begin{figure}[htbp]
\begin{center}
\unitlength 0.1in
\begin{picture}( 45.0000, 20.1000)(  5.0000,-25.0500)
%
\special{pn 8}%
\special{ar 1500 1500 1000 1000  1.5707963 4.7123890}%
%
\special{pn 8}%
\special{pa 1500 500}%
\special{pa 4000 500}%
\special{fp}%
%
\special{pn 8}%
\special{pa 1500 2500}%
\special{pa 4000 2500}%
\special{fp}%
%
\special{pn 8}%
\special{ar 4000 1500 1000 1000  4.7123890 6.2831853}%
\special{ar 4000 1500 1000 1000  0.0000000 1.5707963}%
%
\special{pn 8}%
\special{pa 2750 1500}%
\special{pa 4000 1500}%
\special{fp}%
%
\special{pn 8}%
\special{ar 5000 1500 1390 1500  2.4087776 3.8744078}%
%
\special{pn 8}%
\special{ar 5000 1500 1760 1500  2.4075040 3.8756813}%
\put(26.7000,-15.6000){\makebox(0,0)[rt]{$0$}}%
\put(18.3000,-20.9000){\makebox(0,0){$U$}}%
\put(30.0000,-15.2000){\makebox(0,0)[lt]{$\gamma$}}%
\put(34.5000,-7.7000){\makebox(0,0)[rb]{$c_1$}}%
\put(38.3000,-8.4000){\makebox(0,0)[lt]{$c_2$}}%
\put(40.8000,-15.6000){\makebox(0,0)[lt]{$x$}}%
%
\special{pn 8}%
\special{sh 0.600}%
\special{ar 4010 1500 28 28  0.0000000 6.2831853}%
%
\special{pn 8}%
\special{sh 0.600}%
\special{ar 2750 1500 28 28  0.0000000 6.2831853}%
\end{picture}%
\caption{}
\label{fig2}
\end{center}
\end{figure}
Then the distance from a point in $\gamma$ to $\Sigma\setminus U$
is a continuous function on $\gamma$, and thus has a positive
minimum. This contradicts the assumption that $\xi$ is
prime.
\qed

\bigskip
A positive valued continuous function $\rho$ on $U$ is called {\em
admissible} if 
$$
\int_U\rho^2dvol<\infty.$$
Given a subset $c$ in $U$, $\rho$-diam$(c)$ denotes the diameter
of $c$ w.\ r.\ t.\ the Riemannian metric $\rho^2 g$. (Function
theorists often denote the same metric by $\rho\abs{dz}$.)
An end $\xi$ is called {\em conformal} if for any admissible function
$\rho$ there is a chain $\{c_i\}$ representing $\xi$
such that $\rho$-diam$(c_i)\to 0$. 

If $\phi:U\to V$ is a
conformal equivalence and if $\rho:V\to(0,\infty)$ is admissible,
then a function $\sigma:U\to(0,\infty)$ defined by 
$\sigma(z)=\rho(\phi(z))\abs{\phi'(z)}$ is admissible, and for 
$c\subset V$, we have $\rho$-diam$(c)=\sigma$-diam$(\phi^{-1}(c))$.
This shows that $\phi$ induces a bijection between the sets
of the conformal ends of the two hyperbolic domains.

\begin{lemma} \label{1}
An end $\xi$ is prime if and only if it is conformal.
\end{lemma}

\bf Proof. \rm First of all assuming that $\xi$ is a prime end which
is represented by a topological chain $\{c_i\}$, we shall show
that $\xi$ is a conformal end. By passing to
a subsequence one may further assume that $\overline c_i$ converges to a point
$x_0$. Since $x_0$ belongs to at most one $\overline c_i$, one may also
assume that $x_0\not\in\overline c_i$ for any $i$.
Take a polar coordinates $(r,\theta)$ around $x_0$.
Let $\rho$ be an arbitrary admissible function on $U$, extended to the
whole
$\Sigma$
by letting  $\rho=0$ outside $U$. Then by the Schwarz inequality
$$
(\int_0^\epsilon\int_0^{2\pi} \rho(r,\theta)rd\theta dr)^2
\leq\pi\epsilon^2\cdot\int_{r\leq\epsilon}\rho^2dvol.
$$
Since $\rho$ is admissible,
 $\int_{r\leq\epsilon}\rho^2dvol\to0$ as $\epsilon\to 0$, 
and we have
$$
\frac{1}{\epsilon}\int_0^\epsilon\int_0^{2\pi} \rho(r,\theta)rd\theta dr
\to 0\ \ (\epsilon\to0).
$$
Therefore
we can find
a sequence $\epsilon_k\downarrow 0$ such that 
$$\int_0^{2\pi} \rho(\epsilon_k,\theta)\epsilon_k d\theta\to 
0\ \ (k\to\infty).
$$ 
Notice that the LHS above coincides with the $\rho$-length
of the union of arcs $\{r=\epsilon_k\}\cap U$.

Now from the sequences $\{c_i\}$ and
$\{\epsilon_k\}$, let us construct subsequences $\{c'_i\}$ and 
$\{\epsilon'_k\}$ by the following fashion. 
See Figure \ref{fig3}. 
\begin{figure}[htbp]
\begin{center}
\unitlength 0.1in
\begin{picture}( 41.7000, 34.0000)(  8.3000,-40.0000)
%
\special{pn 8}%
\special{ar 3000 4000 2000 2000  3.1415927 6.2831853}%
%
\special{pn 8}%
\special{ar 2100 3700 500 500  3.1415927 6.2831853}%
%
\special{pn 8}%
\special{pa 2100 3700}%
\special{pa 2096 3668}%
\special{pa 2092 3636}%
\special{pa 2090 3604}%
\special{pa 2088 3572}%
\special{pa 2090 3540}%
\special{pa 2096 3510}%
\special{pa 2104 3478}%
\special{pa 2114 3448}%
\special{pa 2126 3418}%
\special{pa 2142 3390}%
\special{pa 2160 3362}%
\special{pa 2180 3336}%
\special{pa 2200 3310}%
\special{pa 2222 3286}%
\special{pa 2246 3266}%
\special{pa 2270 3246}%
\special{pa 2296 3226}%
\special{pa 2322 3210}%
\special{pa 2348 3196}%
\special{pa 2374 3182}%
\special{pa 2404 3168}%
\special{pa 2434 3156}%
\special{pa 2468 3142}%
\special{pa 2504 3128}%
\special{pa 2540 3112}%
\special{pa 2578 3098}%
\special{pa 2612 3082}%
\special{pa 2644 3064}%
\special{pa 2670 3048}%
\special{pa 2688 3032}%
\special{pa 2700 3016}%
\special{pa 2700 3000}%
\special{pa 2690 2984}%
\special{pa 2672 2968}%
\special{pa 2648 2954}%
\special{pa 2618 2940}%
\special{pa 2582 2926}%
\special{pa 2544 2914}%
\special{pa 2502 2902}%
\special{pa 2462 2892}%
\special{pa 2420 2882}%
\special{pa 2382 2874}%
\special{pa 2346 2868}%
\special{pa 2312 2862}%
\special{pa 2280 2858}%
\special{pa 2248 2854}%
\special{pa 2218 2850}%
\special{pa 2190 2846}%
\special{pa 2162 2844}%
\special{pa 2134 2840}%
\special{pa 2106 2838}%
\special{pa 2078 2834}%
\special{pa 2052 2832}%
\special{pa 2024 2828}%
\special{pa 1996 2824}%
\special{pa 1966 2820}%
\special{pa 1938 2816}%
\special{pa 1908 2812}%
\special{pa 1878 2808}%
\special{pa 1846 2802}%
\special{pa 1812 2798}%
\special{pa 1778 2792}%
\special{pa 1742 2786}%
\special{pa 1704 2780}%
\special{pa 1664 2774}%
\special{pa 1622 2768}%
\special{pa 1580 2760}%
\special{pa 1536 2752}%
\special{pa 1492 2744}%
\special{pa 1448 2736}%
\special{pa 1406 2728}%
\special{pa 1364 2718}%
\special{pa 1328 2708}%
\special{pa 1294 2696}%
\special{pa 1264 2686}%
\special{pa 1240 2674}%
\special{pa 1222 2660}%
\special{pa 1212 2648}%
\special{pa 1208 2632}%
\special{pa 1212 2618}%
\special{pa 1224 2602}%
\special{pa 1246 2586}%
\special{pa 1272 2568}%
\special{pa 1302 2552}%
\special{pa 1336 2536}%
\special{pa 1372 2520}%
\special{pa 1408 2506}%
\special{pa 1442 2492}%
\special{pa 1478 2480}%
\special{pa 1510 2470}%
\special{pa 1542 2460}%
\special{pa 1574 2450}%
\special{pa 1604 2442}%
\special{pa 1634 2434}%
\special{pa 1664 2428}%
\special{pa 1692 2422}%
\special{pa 1722 2414}%
\special{pa 1752 2408}%
\special{pa 1780 2402}%
\special{pa 1810 2396}%
\special{pa 1840 2390}%
\special{pa 1870 2384}%
\special{pa 1900 2378}%
\special{pa 1932 2370}%
\special{pa 1962 2364}%
\special{pa 1994 2358}%
\special{pa 2026 2350}%
\special{pa 2058 2344}%
\special{pa 2088 2336}%
\special{pa 2120 2330}%
\special{pa 2152 2324}%
\special{pa 2184 2316}%
\special{pa 2216 2310}%
\special{pa 2248 2304}%
\special{pa 2280 2296}%
\special{pa 2312 2290}%
\special{pa 2344 2284}%
\special{pa 2374 2278}%
\special{pa 2406 2272}%
\special{pa 2438 2266}%
\special{pa 2468 2258}%
\special{pa 2500 2252}%
\special{pa 2530 2246}%
\special{pa 2560 2238}%
\special{pa 2592 2230}%
\special{pa 2622 2220}%
\special{pa 2652 2210}%
\special{pa 2682 2198}%
\special{pa 2712 2186}%
\special{pa 2742 2172}%
\special{pa 2770 2158}%
\special{pa 2800 2140}%
\special{pa 2828 2122}%
\special{pa 2858 2102}%
\special{pa 2886 2080}%
\special{pa 2912 2056}%
\special{pa 2938 2032}%
\special{pa 2962 2006}%
\special{pa 2982 1980}%
\special{pa 3000 1954}%
\special{pa 3014 1928}%
\special{pa 3026 1902}%
\special{pa 3032 1876}%
\special{pa 3032 1850}%
\special{pa 3030 1826}%
\special{pa 3022 1802}%
\special{pa 3010 1778}%
\special{pa 2996 1756}%
\special{pa 2978 1734}%
\special{pa 2956 1712}%
\special{pa 2932 1692}%
\special{pa 2904 1672}%
\special{pa 2874 1652}%
\special{pa 2844 1634}%
\special{pa 2810 1616}%
\special{pa 2776 1598}%
\special{pa 2740 1582}%
\special{pa 2704 1566}%
\special{pa 2666 1550}%
\special{pa 2628 1536}%
\special{pa 2590 1522}%
\special{pa 2552 1508}%
\special{pa 2514 1496}%
\special{pa 2478 1484}%
\special{pa 2442 1472}%
\special{pa 2406 1462}%
\special{pa 2372 1452}%
\special{pa 2338 1442}%
\special{pa 2306 1434}%
\special{pa 2272 1424}%
\special{pa 2240 1416}%
\special{pa 2208 1410}%
\special{pa 2178 1402}%
\special{pa 2148 1394}%
\special{pa 2116 1388}%
\special{pa 2086 1382}%
\special{pa 2056 1374}%
\special{pa 2028 1368}%
\special{pa 1998 1362}%
\special{pa 1968 1356}%
\special{pa 1940 1348}%
\special{pa 1910 1342}%
\special{pa 1882 1336}%
\special{pa 1852 1328}%
\special{pa 1824 1322}%
\special{pa 1794 1314}%
\special{pa 1766 1306}%
\special{pa 1736 1298}%
\special{pa 1706 1290}%
\special{pa 1676 1280}%
\special{pa 1646 1270}%
\special{pa 1616 1260}%
\special{pa 1586 1250}%
\special{pa 1556 1238}%
\special{pa 1526 1228}%
\special{pa 1496 1214}%
\special{pa 1466 1202}%
\special{pa 1436 1188}%
\special{pa 1406 1174}%
\special{pa 1378 1160}%
\special{pa 1350 1144}%
\special{pa 1322 1128}%
\special{pa 1294 1112}%
\special{pa 1266 1094}%
\special{pa 1240 1076}%
\special{pa 1214 1056}%
\special{pa 1190 1036}%
\special{pa 1166 1016}%
\special{pa 1142 994}%
\special{pa 1120 972}%
\special{pa 1096 950}%
\special{pa 1076 926}%
\special{pa 1054 902}%
\special{pa 1034 878}%
\special{pa 1012 852}%
\special{pa 992 828}%
\special{pa 972 802}%
\special{pa 954 776}%
\special{pa 934 750}%
\special{pa 916 724}%
\special{pa 898 698}%
\special{pa 878 670}%
\special{pa 860 644}%
\special{pa 842 618}%
\special{pa 830 600}%
\special{sp}%
%
\special{pn 8}%
\special{pa 2100 3700}%
\special{pa 2112 3670}%
\special{pa 2124 3640}%
\special{pa 2136 3610}%
\special{pa 2148 3582}%
\special{pa 2164 3554}%
\special{pa 2180 3526}%
\special{pa 2198 3498}%
\special{pa 2218 3474}%
\special{pa 2240 3450}%
\special{pa 2264 3428}%
\special{pa 2290 3408}%
\special{pa 2316 3392}%
\special{pa 2344 3376}%
\special{pa 2374 3364}%
\special{pa 2404 3354}%
\special{pa 2434 3344}%
\special{pa 2466 3336}%
\special{pa 2496 3328}%
\special{pa 2528 3322}%
\special{pa 2560 3314}%
\special{pa 2592 3304}%
\special{pa 2624 3296}%
\special{pa 2654 3286}%
\special{pa 2684 3274}%
\special{pa 2714 3264}%
\special{pa 2744 3252}%
\special{pa 2774 3240}%
\special{pa 2802 3228}%
\special{pa 2832 3216}%
\special{pa 2860 3204}%
\special{pa 2888 3190}%
\special{pa 2914 3176}%
\special{pa 2942 3160}%
\special{pa 2970 3142}%
\special{pa 2998 3122}%
\special{pa 3028 3100}%
\special{pa 3054 3074}%
\special{pa 3074 3048}%
\special{pa 3088 3024}%
\special{pa 3090 2998}%
\special{pa 3082 2976}%
\special{pa 3066 2954}%
\special{pa 3040 2932}%
\special{pa 3010 2914}%
\special{pa 2976 2896}%
\special{pa 2938 2880}%
\special{pa 2902 2864}%
\special{pa 2866 2852}%
\special{pa 2832 2840}%
\special{pa 2800 2828}%
\special{pa 2770 2818}%
\special{pa 2740 2810}%
\special{pa 2712 2802}%
\special{pa 2684 2796}%
\special{pa 2656 2790}%
\special{pa 2628 2784}%
\special{pa 2600 2778}%
\special{pa 2572 2774}%
\special{pa 2542 2768}%
\special{pa 2512 2762}%
\special{pa 2480 2758}%
\special{pa 2448 2750}%
\special{pa 2412 2742}%
\special{pa 2376 2734}%
\special{pa 2336 2724}%
\special{pa 2296 2712}%
\special{pa 2256 2700}%
\special{pa 2220 2688}%
\special{pa 2188 2674}%
\special{pa 2162 2658}%
\special{pa 2146 2644}%
\special{pa 2140 2630}%
\special{pa 2146 2616}%
\special{pa 2164 2600}%
\special{pa 2188 2586}%
\special{pa 2222 2574}%
\special{pa 2258 2560}%
\special{pa 2298 2548}%
\special{pa 2340 2536}%
\special{pa 2382 2526}%
\special{pa 2420 2514}%
\special{pa 2458 2506}%
\special{pa 2496 2496}%
\special{pa 2532 2488}%
\special{pa 2566 2480}%
\special{pa 2598 2472}%
\special{pa 2632 2466}%
\special{pa 2662 2458}%
\special{pa 2692 2452}%
\special{pa 2722 2446}%
\special{pa 2752 2438}%
\special{pa 2780 2432}%
\special{pa 2808 2426}%
\special{pa 2836 2418}%
\special{pa 2864 2410}%
\special{pa 2890 2402}%
\special{pa 2918 2394}%
\special{pa 2944 2386}%
\special{pa 2970 2376}%
\special{pa 2998 2366}%
\special{pa 3024 2354}%
\special{pa 3052 2342}%
\special{pa 3078 2330}%
\special{pa 3106 2318}%
\special{pa 3134 2304}%
\special{pa 3162 2288}%
\special{pa 3190 2272}%
\special{pa 3216 2256}%
\special{pa 3244 2240}%
\special{pa 3272 2222}%
\special{pa 3298 2204}%
\special{pa 3326 2186}%
\special{pa 3354 2166}%
\special{pa 3380 2146}%
\special{pa 3406 2126}%
\special{pa 3432 2106}%
\special{pa 3458 2084}%
\special{pa 3484 2064}%
\special{pa 3508 2040}%
\special{pa 3534 2018}%
\special{pa 3558 1996}%
\special{pa 3582 1972}%
\special{pa 3604 1948}%
\special{pa 3628 1924}%
\special{pa 3650 1900}%
\special{pa 3670 1876}%
\special{pa 3692 1850}%
\special{pa 3712 1826}%
\special{pa 3732 1800}%
\special{pa 3752 1774}%
\special{pa 3770 1746}%
\special{pa 3788 1720}%
\special{pa 3804 1692}%
\special{pa 3820 1666}%
\special{pa 3836 1638}%
\special{pa 3852 1610}%
\special{pa 3866 1580}%
\special{pa 3878 1552}%
\special{pa 3890 1522}%
\special{pa 3902 1492}%
\special{pa 3914 1462}%
\special{pa 3922 1432}%
\special{pa 3932 1402}%
\special{pa 3940 1372}%
\special{pa 3948 1340}%
\special{pa 3954 1310}%
\special{pa 3960 1278}%
\special{pa 3964 1246}%
\special{pa 3968 1214}%
\special{pa 3972 1182}%
\special{pa 3974 1150}%
\special{pa 3978 1118}%
\special{pa 3978 1086}%
\special{pa 3980 1054}%
\special{pa 3980 1020}%
\special{pa 3982 988}%
\special{pa 3980 956}%
\special{pa 3980 924}%
\special{pa 3980 892}%
\special{pa 3978 860}%
\special{pa 3976 828}%
\special{pa 3974 796}%
\special{pa 3972 764}%
\special{pa 3970 732}%
\special{pa 3968 700}%
\special{pa 3966 668}%
\special{pa 3964 636}%
\special{pa 3960 604}%
\special{pa 3960 600}%
\special{sp}%
%
\special{pn 13}%
\special{pa 2180 1400}%
\special{pa 3940 1400}%
\special{fp}%
%
\special{pn 20}%
\special{ar 3000 4000 2000 2000  4.6970645 4.9553206}%
%
\special{pn 20}%
\special{pa 2300 2870}%
\special{pa 2530 2770}%
\special{fp}%
%
\special{pn 20}%
\special{ar 2100 3700 500 500  5.0611600 5.4777898}%
%
\special{pn 13}%
\special{ar 2100 3700 250 250  4.7531827 5.1671083}%
%
\special{pn 13}%
\special{sh 1}%
\special{ar 2100 3700 10 10 0  6.28318530717959E+0000}%
\put(25.8000,-8.3000){\makebox(0,0){$U$}}%
\put(30.4000,-13.6000){\makebox(0,0)[lb]{$c_1'$}}%
\put(31.0000,-19.9000){\makebox(0,0)[lb]{$c_1''$}}%
\put(25.7000,-27.5000){\makebox(0,0)[lb]{$c_2'$}}%
\put(23.8000,-32.9000){\makebox(0,0)[lb]{$c_2''$}}%
\put(20.5000,-34.4000){\makebox(0,0)[rt]{$c_3'$}}%
\put(26.0000,-35.9000){\makebox(0,0)[lb]{$r=\varepsilon_2'$}}%
\put(44.0000,-25.9000){\makebox(0,0)[lb]{$r=\varepsilon_1'$}}%
\put(21.2000,-37.2000){\makebox(0,0)[lt]{$x_0$}}%
\end{picture}%
\caption{}
\label{fig3}
\end{center}
\end{figure}
First define
$c'_1=c_1$ and choose $\epsilon'_1$ to be any $\epsilon_k$ from
the sequence such that $\overline c'_1\cap\{r\leq\epsilon'_1\}=\emptyset$.
Then choose $c'_2$ to be any $c_i$ from the sequence such that
$\overline c'_2\subset\{r<\epsilon'_1\}$. Next choose $\epsilon'_2$
such that $\overline c'_2\cap\{r\leq\epsilon'_2\}=\emptyset$,
$c'_3$ such that $\overline c'_3\subset\{r<\epsilon'_2\}$,
and so forth. 

Then there is a connected component $c''_i$ 
of $\{r=\epsilon'_i\}\cap U$ which separates the cross cut $c'_{i+1}$
from $c'_i$. To see this, construct a graph $\Gamma$; the vertices are
connected components of $U\setminus\{r=\epsilon'_i\}$ and the edges
connected components of $U\cap\{r=\epsilon'_i\}$. See figure \ref{fig4}.
\begin{figure}[htbp]
\begin{center}
\unitlength 0.1in
%
\caption{}
\label{fig4}
\end{center}
\end{figure}
By transversality
argument any two distinct vertices can be joined by a finite edge path.
Actually $\Gamma$ is a tree, since $U$ is simply connected and any edge
corresponds to a cross cut of $U$. Thus there is a unique shortest
edge path joining the two vertices corresponding to the components,
one containing $c'_i$, the other $c'_{i+1}$.
The component $c''_i$
of $U\cap\{r=\epsilon'_i\}$ corresponding to any edge of $\sigma$ 
separates $c'_{i+1}$ from $c'_i$.
Clearly the chains $\{c'_i\}$ and $\{c''_i\}$
are equivalent and the latter satisfies $\rho$-diam$(c''_i)\to 0$,
showing that $\xi$ is conformal.

Next assume that $\xi$ is conformal.
First of all if we choose an admissible
function $\rho_0$ which is constantly equal to 1 on $U$, we can find
a chain $\{c_i\}$ such that ${\rm diam}(c_i)\to 0$ as $i\to\infty$.
Passing to a subsequence if necessary, one may assume $c_i\to x_0$.
Again let $(r,\theta)$ be the polar coordinates around $x_0$.
Define a function $\rho$ by 
$$
\rho(r,\theta)=-\frac{1}{r\log r}
$$
if $r\leq 1/2$ and equal to $2/\log 2$ otherwise.
Computation shows that the restriction of $\rho$ to $U$ is admissible.
Now for any small $\epsilon>\delta$, the $\rho$-distance 
of the $\epsilon$-circle and the $\delta$-circle is given by
$$-\int_\delta^\epsilon\frac{dr}{r\log r}=\log(\log\delta/\log\epsilon),
$$ which diverges to $\infty$ if we fix $\epsilon$ and let $\delta\to0$.
Let $c'_i$ be a chain representing $\xi$ such that
$\rho$-diam$(c'_i)\to 0$. Since $\rho$ is bigger than a constant
multiple of $\rho_0$, this implies also that diam$(c'_i)\to 0$.

First consider the case where $c'_i$ converges to $x_0$ (passing
to a subsequence). See Figure \ref{fig5}.
\begin{figure}[htbp]
\begin{center}
\unitlength 0.1in
\begin{picture}( 44.7000, 21.9600)(  5.0000,-25.1600)
%
\special{pn 8}%
\special{sh 0.600}%
\special{ar 1500 1500 28 28  0.0000000 6.2831853}%
%
\special{pn 8}%
\special{ar 1500 1500 1000 1000  1.5707963 4.7123890}%
%
\special{pn 8}%
\special{sh 0.600}%
\special{ar 4900 1500 28 28  0.0000000 6.2831853}%
%
\special{pn 8}%
\special{pa 1510 500}%
\special{pa 1542 502}%
\special{pa 1574 504}%
\special{pa 1606 508}%
\special{pa 1638 510}%
\special{pa 1670 514}%
\special{pa 1702 518}%
\special{pa 1734 524}%
\special{pa 1764 532}%
\special{pa 1796 540}%
\special{pa 1826 550}%
\special{pa 1856 560}%
\special{pa 1886 572}%
\special{pa 1916 586}%
\special{pa 1944 598}%
\special{pa 1974 612}%
\special{pa 2002 626}%
\special{pa 2032 640}%
\special{pa 2060 654}%
\special{pa 2088 668}%
\special{pa 2118 684}%
\special{pa 2146 698}%
\special{pa 2174 714}%
\special{pa 2202 728}%
\special{pa 2230 744}%
\special{pa 2258 760}%
\special{pa 2286 774}%
\special{pa 2314 790}%
\special{pa 2342 806}%
\special{pa 2370 822}%
\special{pa 2398 836}%
\special{pa 2426 852}%
\special{pa 2454 868}%
\special{pa 2482 884}%
\special{pa 2510 898}%
\special{pa 2538 914}%
\special{pa 2566 928}%
\special{pa 2596 942}%
\special{pa 2624 956}%
\special{pa 2654 970}%
\special{pa 2682 984}%
\special{pa 2712 996}%
\special{pa 2740 1010}%
\special{pa 2770 1020}%
\special{pa 2800 1032}%
\special{pa 2830 1044}%
\special{pa 2860 1054}%
\special{pa 2890 1064}%
\special{pa 2922 1074}%
\special{pa 2952 1084}%
\special{pa 2982 1092}%
\special{pa 3014 1102}%
\special{pa 3044 1110}%
\special{pa 3074 1118}%
\special{pa 3106 1126}%
\special{pa 3138 1134}%
\special{pa 3168 1142}%
\special{pa 3200 1150}%
\special{pa 3230 1156}%
\special{pa 3262 1164}%
\special{pa 3294 1170}%
\special{pa 3324 1178}%
\special{pa 3356 1186}%
\special{pa 3388 1192}%
\special{pa 3418 1200}%
\special{pa 3450 1206}%
\special{pa 3482 1214}%
\special{pa 3512 1220}%
\special{pa 3544 1228}%
\special{pa 3576 1234}%
\special{pa 3606 1242}%
\special{pa 3638 1248}%
\special{pa 3670 1256}%
\special{pa 3700 1262}%
\special{pa 3732 1268}%
\special{pa 3762 1276}%
\special{pa 3794 1282}%
\special{pa 3826 1290}%
\special{pa 3856 1296}%
\special{pa 3888 1302}%
\special{pa 3918 1310}%
\special{pa 3950 1316}%
\special{pa 3982 1322}%
\special{pa 4012 1328}%
\special{pa 4044 1336}%
\special{pa 4076 1342}%
\special{pa 4106 1348}%
\special{pa 4138 1354}%
\special{pa 4170 1360}%
\special{pa 4200 1366}%
\special{pa 4232 1372}%
\special{pa 4264 1378}%
\special{pa 4296 1384}%
\special{pa 4326 1390}%
\special{pa 4358 1396}%
\special{pa 4390 1402}%
\special{pa 4420 1408}%
\special{pa 4452 1414}%
\special{pa 4484 1420}%
\special{pa 4516 1426}%
\special{pa 4546 1432}%
\special{pa 4578 1438}%
\special{pa 4590 1440}%
\special{sp}%
%
\special{pn 8}%
\special{pa 1500 2500}%
\special{pa 1532 2498}%
\special{pa 1564 2496}%
\special{pa 1596 2494}%
\special{pa 1628 2490}%
\special{pa 1660 2488}%
\special{pa 1692 2482}%
\special{pa 1724 2476}%
\special{pa 1754 2470}%
\special{pa 1786 2460}%
\special{pa 1816 2452}%
\special{pa 1846 2440}%
\special{pa 1876 2428}%
\special{pa 1906 2416}%
\special{pa 1934 2402}%
\special{pa 1964 2390}%
\special{pa 1992 2376}%
\special{pa 2022 2360}%
\special{pa 2050 2346}%
\special{pa 2078 2332}%
\special{pa 2108 2318}%
\special{pa 2136 2302}%
\special{pa 2164 2288}%
\special{pa 2192 2272}%
\special{pa 2220 2256}%
\special{pa 2248 2242}%
\special{pa 2276 2226}%
\special{pa 2304 2210}%
\special{pa 2332 2196}%
\special{pa 2360 2180}%
\special{pa 2388 2164}%
\special{pa 2416 2148}%
\special{pa 2444 2134}%
\special{pa 2472 2118}%
\special{pa 2500 2102}%
\special{pa 2528 2088}%
\special{pa 2556 2074}%
\special{pa 2586 2058}%
\special{pa 2614 2044}%
\special{pa 2644 2030}%
\special{pa 2672 2018}%
\special{pa 2702 2004}%
\special{pa 2730 1992}%
\special{pa 2760 1980}%
\special{pa 2790 1968}%
\special{pa 2820 1958}%
\special{pa 2850 1948}%
\special{pa 2880 1936}%
\special{pa 2912 1928}%
\special{pa 2942 1918}%
\special{pa 2972 1908}%
\special{pa 3004 1900}%
\special{pa 3034 1892}%
\special{pa 3064 1882}%
\special{pa 3096 1874}%
\special{pa 3128 1866}%
\special{pa 3158 1860}%
\special{pa 3190 1852}%
\special{pa 3220 1844}%
\special{pa 3252 1838}%
\special{pa 3284 1830}%
\special{pa 3314 1822}%
\special{pa 3346 1816}%
\special{pa 3378 1808}%
\special{pa 3408 1802}%
\special{pa 3440 1794}%
\special{pa 3472 1788}%
\special{pa 3502 1780}%
\special{pa 3534 1774}%
\special{pa 3566 1766}%
\special{pa 3596 1760}%
\special{pa 3628 1752}%
\special{pa 3660 1746}%
\special{pa 3690 1740}%
\special{pa 3722 1732}%
\special{pa 3752 1726}%
\special{pa 3784 1718}%
\special{pa 3816 1712}%
\special{pa 3846 1704}%
\special{pa 3878 1698}%
\special{pa 3908 1692}%
\special{pa 3940 1686}%
\special{pa 3972 1678}%
\special{pa 4002 1672}%
\special{pa 4034 1666}%
\special{pa 4066 1660}%
\special{pa 4096 1654}%
\special{pa 4128 1646}%
\special{pa 4160 1640}%
\special{pa 4190 1634}%
\special{pa 4222 1628}%
\special{pa 4254 1622}%
\special{pa 4286 1616}%
\special{pa 4316 1610}%
\special{pa 4348 1604}%
\special{pa 4380 1598}%
\special{pa 4410 1592}%
\special{pa 4442 1586}%
\special{pa 4474 1580}%
\special{pa 4506 1574}%
\special{pa 4536 1568}%
\special{pa 4568 1562}%
\special{pa 4580 1560}%
\special{sp}%
%
\special{pn 8}%
\special{ar 4900 1500 2650 2650  2.7481420 3.5350433}%
%
\special{pn 8}%
\special{ar 4900 1500 1164 1164  2.9578481 3.3406746}%
%
\special{pn 8}%
\special{ar 4900 1500 1848 1848  2.9410292 3.3549613}%
\put(15.2000,-15.3000){\makebox(0,0)[lt]{$0$}}%
\put(10.4000,-15.0000){\makebox(0,0){$U$}}%
\put(24.7000,-4.9000){\makebox(0,0)[lb]{$r=1/2$}}%
\put(30.1000,-15.5000){\makebox(0,0)[rb]{$\overline{c}_i'$}}%
\put(38.0000,-12.4000){\makebox(0,0)[lb]{$\overline{c}_{i+1}'$}}%
\put(49.7000,-15.1000){\makebox(0,0)[lt]{$x_0$}}%
\end{picture}%
\caption{$U$ with the metric $\rho g^2$. $(\rho$-$\mathrm{diam}\,c_i'\to 0$)}
\label{fig5}
\end{center}
\end{figure}
The above computation shows that for $i$ big enough
$\overline c'_i$ is a compact subset of $\{0<r<1/2\}$ and  we
can take a subsequence such that the closures  $\overline c'_i$ are
mutually disjoint. Thus we obtain a topological chain representing
$\xi$.

In the remaining case, we may assume that $c'_i$ converges to
a point $x_1$ distinct from $x_0$. See Figure \ref{fig6}. 
\begin{figure}[htbp]
\begin{center}
\unitlength 0.1in
\begin{picture}( 49.2500, 33.6900)(  4.1000,-35.0500)
%
\special{pn 8}%
\special{ar 2500 3000 500 500  3.1415927 6.2831853}%
\special{ar 2500 3000 500 500  0.0000000 1.5707963}%
%
\special{pn 8}%
\special{sh 0.600}%
\special{ar 2500 3000 28 28  0.0000000 6.2831853}%
%
\special{pn 8}%
\special{sh 0.600}%
\special{ar 4800 3000 28 28  0.0000000 6.2831853}%
%
\special{pn 8}%
\special{pa 2500 2700}%
\special{pa 5000 2700}%
\special{fp}%
%
\special{pn 8}%
\special{pa 2500 2300}%
\special{pa 5000 2300}%
\special{fp}%
%
\special{pn 8}%
\special{pa 2500 2400}%
\special{pa 5000 2400}%
\special{fp}%
%
\special{pn 8}%
\special{pa 2500 2600}%
\special{pa 5000 2600}%
\special{fp}%
%
\special{pn 8}%
\special{pa 5000 2900}%
\special{pa 4500 2900}%
\special{fp}%
%
\special{pn 8}%
\special{pa 5000 2800}%
\special{pa 4500 2800}%
\special{fp}%
%
\special{pn 8}%
\special{pa 2500 2200}%
\special{pa 5000 2200}%
\special{fp}%
%
\special{pn 8}%
\special{pa 2500 2100}%
\special{pa 5000 2100}%
\special{fp}%
%
\special{pn 8}%
\special{ar 5000 2750 50 50  4.7123890 6.2831853}%
\special{ar 5000 2750 50 50  0.0000000 1.5707963}%
%
\special{pn 8}%
\special{ar 5000 2750 150 150  4.7123890 6.2831853}%
\special{ar 5000 2750 150 150  0.0000000 1.5707963}%
%
\special{pn 8}%
\special{ar 2500 2500 100 100  1.5707963 4.7123890}%
%
\special{pn 8}%
\special{ar 2500 2500 200 200  1.5707963 4.7123890}%
%
\special{pn 8}%
\special{ar 5000 2250 50 50  4.7123890 6.2831853}%
\special{ar 5000 2250 50 50  0.0000000 1.5707963}%
%
\special{pn 8}%
\special{ar 5000 2250 150 150  4.7123890 6.2831853}%
\special{ar 5000 2250 150 150  0.0000000 1.5707963}%
%
\special{pn 8}%
\special{ar 2500 2000 100 100  1.5707963 4.7123890}%
%
\special{pn 8}%
\special{ar 2450 430 294 294  1.7506498 4.8821673}%
%
\special{pn 8}%
\special{pa 2400 720}%
\special{pa 4900 1160}%
\special{fp}%
%
\special{pn 8}%
\special{ar 4850 1250 104 104  4.7123890 6.2831853}%
\special{ar 4850 1250 104 104  0.0000000 1.6295521}%
%
\special{pn 8}%
\special{sh 0.600}%
\special{ar 2510 430 28 28  0.0000000 6.2831853}%
%
\special{pn 8}%
\special{ar 2500 1950 250 250  1.5707963 4.4768440}%
%
\special{pn 8}%
\special{pa 4880 1350}%
\special{pa 2420 1710}%
\special{fp}%
%
\special{pn 8}%
\special{ar 4700 970 502 502  0.8139618 0.8276266}%
%
\special{pn 8}%
\special{ar 4820 1070 516 516  4.9332178 6.2831853}%
\special{ar 4820 1070 516 516  0.0000000 1.2551807}%
%
\special{pn 8}%
\special{pa 2500 1890}%
\special{pa 4970 1560}%
\special{fp}%
%
\special{pn 13}%
\special{pa 4800 2700}%
\special{pa 4800 2600}%
\special{fp}%
%
\special{pn 13}%
\special{pa 4800 2200}%
\special{pa 4800 2100}%
\special{fp}%
%
\special{pn 13}%
\special{pa 4800 1580}%
\special{pa 4800 1360}%
\special{fp}%
\put(48.2000,-16.2000){\makebox(0,0)[lt]{$c_1'$}}%
\put(48.2000,-20.8000){\makebox(0,0)[lb]{$c_2'$}}%
\put(48.2000,-25.7000){\makebox(0,0)[lb]{$c_3'$}}%
\put(48.2000,-30.4000){\makebox(0,0)[lt]{$x_1$}}%
\put(25.4000,-30.3000){\makebox(0,0)[lt]{$x_0$}}%
%
\special{pn 20}%
\special{ar 2500 3000 500 500  4.3128735 4.5307402}%
%
\special{pn 13}%
\special{pa 2700 2400}%
\special{pa 2700 2300}%
\special{fp}%
%
\special{pn 20}%
\special{ar 2500 3000 1000 1000  5.1581013 5.3638422}%
%
\special{pn 13}%
\special{pa 2700 1860}%
\special{pa 2700 1670}%
\special{fp}%
%
\special{pn 20}%
\special{ar 2500 3000 1500 1500  5.1040889 5.3209111}%
%
\special{pn 8}%
\special{ar 2500 3000 1510 1510  3.1415927 6.2831853}%
\special{ar 2500 3000 1510 1510  0.0000000 0.3407757}%
%
\special{pn 8}%
\special{ar 2500 3000 1000 1000  3.1415927 6.2831853}%
\special{ar 2500 3000 1000 1000  0.0000000 0.5266273}%
\put(33.7000,-35.7000){\makebox(0,0){$r=\varepsilon_2$}}%
\put(31.3000,-5.1000){\makebox(0,0)[lt]{$U$}}%
\put(30.1000,-15.6000){\makebox(0,0)[lb]{$c_1''$}}%
\put(26.8000,-16.6000){\makebox(0,0)[rb]{$c_1$}}%
\put(31.3000,-20.0000){\makebox(0,0)[lb]{$c_2''$}}%
%
\special{pn 8}%
\special{pa 3120 1990}%
\special{pa 3070 2180}%
\special{fp}%
\special{sh 1}%
\special{pa 3070 2180}%
\special{pa 3106 2122}%
\special{pa 3084 2128}%
\special{pa 3068 2110}%
\special{pa 3070 2180}%
\special{fp}%
\put(27.8000,-24.2000){\makebox(0,0)[lt]{$c_2$}}%
\put(22.0000,-23.4000){\makebox(0,0)[rb]{$c_3''$}}%
%
\special{pn 8}%
\special{pa 2180 2330}%
\special{pa 2370 2520}%
\special{fp}%
\special{sh 1}%
\special{pa 2370 2520}%
\special{pa 2338 2460}%
\special{pa 2332 2482}%
\special{pa 2310 2488}%
\special{pa 2370 2520}%
\special{fp}%
\put(19.4000,-15.8000){\makebox(0,0)[rb]{$r=\varepsilon_1$}}%
\put(25.4000,-4.7000){\makebox(0,0)[lt]{$0$}}%
%
\special{pn 8}%
\special{pa 2500 140}%
\special{pa 4930 570}%
\special{fp}%
\end{picture}%
\caption{}
\label{fig6}
\end{center}
\end{figure}
We shall still use the
polar coordinates $(r,\theta)$ around $x_0$. Recall that
we have another chain $\{c_i\}$ converging to $x_0$.
The chains $\{c_i\}$ has no particularly good property other than
diam$(c_i)\to0$. In the worst case $x_0$ may belong to any
$\overline c_i$.
However passing to subsequences of $\{c'_i\}$ and $\{c_i\}$ (denoted by
the same letters) and
choosing a sequence of positive numbers $\epsilon_i\downarrow 0$, 
we may assume
the following.
\\
(1) The cross cut $c_i$ is contained in $\{r<\epsilon_i\}$.
\\
(2) All the $c'_i$ is disjoint from $\{r\leq\epsilon_1\}$.
\\
(3) The sequence $c'_1,c_1,c'_2,c_2,\cdots$ forms a chain.

Then there is a component $c''_i$ of 
$\{r=\epsilon_{i}\}\cap U$ which separates $c_i$ from $c'_i$.
The chain $\{c''_i\}$ is the desired topological chain.
\qed

\bigskip
A cross cut $c:\R\to U$  is called {\em extendable} if the limits
$\lim_{t\to-\infty}c(t)$ and $\lim_{t\to\infty}c(t)$ exist.
Then $\overline c$ is either a compact arc or a Jordan curve in $\Sigma$.
A topological chain $\{c_i\}$ is called {\em extendable} if each $c_i$
is extendable.
The proof of the above lemma also shows the following lemma
useful in the sequel.

\begin{lemma} \label{1a}
A prime end is represented by an extendable topological chain.
\end{lemma}

\bigskip
For a hyperbolic domain $U$ of $\Sigma$, denote by $\PP(U)$
the set of prime ends of $U$. The union $\hat U=U\cup \PP(U)$,
topologized in a standard way, is called the {\em Carath\'eodory
compactification} of $U$. Let us explain it in bit more details.
A neighbourhood system in $\hat U$ of a point in $U$ is the same
as a given system in $U$. Choose a point $\xi\in \PP(U)$ represented
by a topological chain $\{c_i\}$. The set of points in the content
$U(c_i)$, together with the prime ends represented by 
topological chains contained in $U(c_i)$ for each $i$ forms
a neighbourhood system of $\xi$.

Lemma \ref{1} shows that a conformal equivalence $\phi:U\to V$
extends to a homeomorphism $\hat\phi:\hat U\to\hat V$. In particular
$\hat U$ is homeomorphic to $\hat \D$
by the natural extension $\hat\phi$ of a Riemann mapping $\phi:U\to\D$,
and for $\D$ it is clear that $\hat \D$ is homeomorphic to the
closed disc $\D\cup\boundary$. Thus $\hat U$ is homeomorphic to
a closed disc for any hyperbolic domain $U$.
On the other hand by the definition of topological chains, 
a homeomorphism $f$ of $U$
which extends to a homeomorphism of the closure $\overline U$ does
extend to a homeomorphism $\hat f$ of the compact disc $\hat U$.
Especially important is the rotation number of the restriction of
$\hat f$ to $\PP(U)$, which is called the {\em Carath\'eodory
rotation number}.

A proper embedding $\gamma:[0,\infty)\to U$ is
called a {\em ray}.  A ray $\gamma$ is said to {\em belong
to} a prime end $\xi$ if $\xi$ is represented by a chain $\{c_i\}$
and for any $i$, there is $t>0$ such that 
$\gamma[t,\infty)\subset U(c_i)$. The ray $\gamma$ is called {\em
extendable}
if the limit $\lim_{t\to\infty}\gamma(t)$, called the {\em end point} of
$\gamma$, exists. The end point of an extendable ray in $U$ belongs
to the frontier ${\rm Fr}(U)$.

A prime end $\xi$ of $U$ is called {\em extendable} if there is
an extendable ray belonging to $\xi$.
Denote by $\EE(U)$ the set of extendable prime ends.

\begin{lemma} \label{4a}
The end points of two extendable rays $\gamma_i$ $(i=1,2)$ 
belonging to the same prime end $\xi$
coincide.
\end{lemma}

{\bf Proof}.
The end point of $\gamma_i$ is the limit point of {\em any} topological chain
representing $\xi$.
\qed

\bigskip

Lemma \ref{4a} enables us to define a natural map 
$\Phi:\EE(U)\to{\rm Fr}(U)$.

\begin{lemma} \label{3}
Any extendable ray belongs to some prime end.
\end{lemma}

{\bf Proof}. Given an extendable ray $\gamma$ with end point 
$x\in{\rm Fr}(U)$, one can construct a topological chain from the
concentric circles centered at $x$, by much the same argument
as in the proof of Lemma \ref{1}.
\qed

\bigskip
The above lemma says that a ray $\gamma$ extendable in $U\subset \Sigma$
is extendable in the closed disc $\hat U$.

By an identification $\hat\phi:\PP(U)\to\boundary$ induced from
a Riemann mapping $\phi:U\to\D$, the Lebesgue
measure on $\boundary$ is transformed to a probability measure
on $\PP(U)$. It depends upon the choice of the Riemann mapping $\phi$, but
its class (called {\em Lebesgue class}) is unique.

\begin{lemma} \label{4}
The set $\EE(U)$ of extendable prime ends is conull w.\ r.\ t.\
the Lebesgue class. Especially $\EE(U)$ is dense in $\PP(U)$.
\end{lemma}

{\bf Proof}. Let $\psi:\D\to U$ be the inverse Riemann
mapping. Then another application of the Schwarz inequality shows
$$
\int_0^{2\pi}\int_{1/2}^1\abs{\psi'(re^{i\theta})}rdrd\theta<\infty.$$
That is, for Lebesgue almost all $\theta_0$, the value
$$2\int_{1/2}^1\abs{\psi'(re^{i\theta_0})}dr
<4\int_{1/2}^1\abs{\psi'(re^{i\theta_0})}rdr<\infty.
$$
Notice that the LHS is the length of the ray
$\psi\{re^{i\theta_0}\mid 1/2\leq r<1\}$.
\qed

\bigskip

\begin{remark}
It is not the case that an extendable prime end always admits a ray
of finite length. See Figure \ref{fig7}.
\begin{figure}[htbp]
\begin{center}
\unitlength 0.1in
\begin{picture}( 30.0000, 15.3900)(  5.0000,-18.4000)
%
\special{pn 8}%
\special{pa 500 440}%
\special{pa 3500 440}%
\special{pa 3500 1840}%
\special{pa 500 1840}%
\special{pa 500 440}%
\special{pa 3500 440}%
\special{fp}%
%
\special{pn 8}%
\special{ar 2000 440 1000 1400  6.1835167 6.2335167}%
\special{ar 2000 440 1000 1400  6.2635167 6.3135167}%
\special{ar 2000 440 1000 1400  6.3435167 6.3935167}%
\special{ar 2000 440 1000 1400  6.4235167 6.4735167}%
\special{ar 2000 440 1000 1400  6.5035167 6.5535167}%
\special{ar 2000 440 1000 1400  6.5835167 6.6335167}%
\special{ar 2000 440 1000 1400  6.6635167 6.7135167}%
\special{ar 2000 440 1000 1400  6.7435167 6.7935167}%
\special{ar 2000 440 1000 1400  6.8235167 6.8735167}%
\special{ar 2000 440 1000 1400  6.9035167 6.9535167}%
\special{ar 2000 440 1000 1400  6.9835167 7.0335167}%
\special{ar 2000 440 1000 1400  7.0635167 7.1135167}%
\special{ar 2000 440 1000 1400  7.1435167 7.1935167}%
\special{ar 2000 440 1000 1400  7.2235167 7.2735167}%
\special{ar 2000 440 1000 1400  7.3035167 7.3535167}%
\special{ar 2000 440 1000 1400  7.3835167 7.4335167}%
\special{ar 2000 440 1000 1400  7.4635167 7.5135167}%
\special{ar 2000 440 1000 1400  7.5435167 7.5935167}%
\special{ar 2000 440 1000 1400  7.6235167 7.6735167}%
\special{ar 2000 440 1000 1400  7.7035167 7.7535167}%
\special{ar 2000 440 1000 1400  7.7835167 7.8335167}%
\special{ar 2000 440 1000 1400  7.8635167 7.9135167}%
\special{ar 2000 440 1000 1400  7.9435167 7.9935167}%
\special{ar 2000 440 1000 1400  8.0235167 8.0735167}%
\special{ar 2000 440 1000 1400  8.1035167 8.1535167}%
\special{ar 2000 440 1000 1400  8.1835167 8.2335167}%
\special{ar 2000 440 1000 1400  8.2635167 8.3135167}%
\special{ar 2000 440 1000 1400  8.3435167 8.3935167}%
\special{ar 2000 440 1000 1400  8.4235167 8.4735167}%
\special{ar 2000 440 1000 1400  8.5035167 8.5535167}%
\special{ar 2000 440 1000 1400  8.5835167 8.6335167}%
\special{ar 2000 440 1000 1400  8.6635167 8.7135167}%
\special{ar 2000 440 1000 1400  8.7435167 8.7935167}%
\special{ar 2000 440 1000 1400  8.8235167 8.8735167}%
\special{ar 2000 440 1000 1400  8.9035167 8.9535167}%
\special{ar 2000 440 1000 1400  8.9835167 9.0335167}%
\special{ar 2000 440 1000 1400  9.0635167 9.1135167}%
\special{ar 2000 440 1000 1400  9.1435167 9.1935167}%
\special{ar 2000 440 1000 1400  9.2235167 9.2735167}%
\special{ar 2000 440 1000 1400  9.3035167 9.3535167}%
\special{ar 2000 440 1000 1400  9.3835167 9.4335167}%
\special{ar 2000 440 1000 1400  9.4635167 9.5135167}%
\put(24.7000,-5.3100){\makebox(0,0){$U$}}%
%
\special{pn 8}%
\special{pa 3500 740}%
\special{pa 1040 740}%
\special{fp}%
%
\special{pn 8}%
\special{pa 500 990}%
\special{pa 2910 990}%
\special{fp}%
%
\special{pn 8}%
\special{pa 3500 1190}%
\special{pa 1170 1190}%
\special{fp}%
%
\special{pn 8}%
\special{pa 500 1340}%
\special{pa 2760 1340}%
\special{fp}%
%
\special{pn 8}%
\special{pa 3500 1440}%
\special{pa 1310 1440}%
\special{fp}%
%
\special{pn 8}%
\special{pa 500 1510}%
\special{pa 2630 1510}%
\special{fp}%
%
\special{pn 8}%
\special{pa 3500 1560}%
\special{pa 1420 1560}%
\special{fp}%
%
\special{pn 8}%
\special{pa 500 1600}%
\special{pa 2550 1600}%
\special{fp}%
%
\special{pn 20}%
\special{pa 2000 590}%
\special{pa 1000 590}%
\special{fp}%
%
\special{pn 20}%
\special{ar 1000 740 150 150  1.5707963 4.7123890}%
%
\special{pn 20}%
\special{pa 1000 890}%
\special{pa 2950 890}%
\special{fp}%
%
\special{pn 20}%
\special{ar 2950 990 100 100  4.7123890 6.2831853}%
\special{ar 2950 990 100 100  0.0000000 1.5707963}%
%
\special{pn 20}%
\special{pa 2950 1090}%
\special{pa 2450 1090}%
\special{fp}%
\special{sh 1}%
\special{pa 2450 1090}%
\special{pa 2518 1110}%
\special{pa 2504 1090}%
\special{pa 2518 1070}%
\special{pa 2450 1090}%
\special{fp}%
\put(8.1000,-7.6000){\makebox(0,0)[rb]{ray}}%
%
\special{pn 8}%
\special{pa 2284 1020}%
\special{pa 2310 1038}%
\special{pa 2330 1064}%
\special{pa 2338 1096}%
\special{pa 2340 1120}%
\special{sp}%
%
\special{pn 8}%
\special{pa 2286 1022}%
\special{pa 2264 1138}%
\special{da 0.070}%
%
\special{pn 8}%
\special{pa 2340 1116}%
\special{pa 2268 1116}%
\special{da 0.070}%
%
\special{pn 8}%
\special{pa 2216 1138}%
\special{pa 2338 1136}%
\special{da 0.070}%
%
\special{pn 8}%
\special{pa 2334 1136}%
\special{pa 2334 1168}%
\special{da 0.070}%
%
\special{pn 8}%
\special{pa 2332 1164}%
\special{pa 2238 1164}%
\special{da 0.070}%
%
\special{pn 8}%
\special{pa 2240 1164}%
\special{pa 2218 1136}%
\special{da 0.070}%
\end{picture}%
\caption{}
\label{fig7}
\end{center}
\end{figure}
\end{remark}

\section{The Cartwright-Littlewood theorem}

Let $f:\Sigma\to \Sigma$ be an orientation preserving homeomorphism which
leaves a hyperbolic domain $U$ in $\Sigma$ invariant. 
Now $f$ induces an
orientation preserving homeomorphism on the Carath\'eodory compactification,
$\hat f:\hat U\to\hat U$. 
The purpose of this section is to give a proof of the following theorem
due to M. L. Cartwright and J. E. Littlewood (\cite{CL}).

\begin{theorem} \label{5}
Let $f$ and $U$ be as above. Assume that there is no fixed point
in ${\rm Fr}(U)$ and that the Carath\'eodory rotation number of $U$ is $0$.
Then the restriction of $\hat f$ to $\PP(U)$ is Morse-Smale, and if
$\xi\in\PP(U)$ is an attractor (resp.\ repellor) of the restriction of
$\hat f$ to $\PP(U)$, then $\xi$ is an attractor (resp.\ repellor)
of the homeomorphism $\hat f$ of $\hat U$.
\end{theorem}

See Figure \ref{fig8}.
\begin{figure}[htbp]
\begin{center}
\unitlength 0.1in
\begin{picture}( 30.5600, 30.5600)(  4.7200,-35.2800)
%
\special{pn 8}%
\special{ar 2000 2000 1500 1500  2.3561945 3.9269908}%
%
\special{pn 8}%
\special{pa 930 950}%
\special{pa 940 940}%
\special{fp}%
\special{sh 1}%
\special{pa 940 940}%
\special{pa 880 976}%
\special{pa 904 978}%
\special{pa 910 1002}%
\special{pa 940 940}%
\special{fp}%
%
\special{pn 8}%
\special{pa 930 3050}%
\special{pa 940 3062}%
\special{fp}%
\special{sh 1}%
\special{pa 940 3062}%
\special{pa 910 2998}%
\special{pa 904 3022}%
\special{pa 880 3026}%
\special{pa 940 3062}%
\special{fp}%
%
\special{pn 8}%
\special{ar 2000 2000 1500 1500  0.0000000 6.2831853}%
%
\special{pn 8}%
\special{ar 2000 2000 1500 1500  5.4977871 6.2831853}%
\special{ar 2000 2000 1500 1500  0.0000000 0.7853982}%
%
\special{pn 8}%
\special{pa 3072 3050}%
\special{pa 3062 3062}%
\special{fp}%
\special{sh 1}%
\special{pa 3062 3062}%
\special{pa 3122 3026}%
\special{pa 3098 3022}%
\special{pa 3092 2998}%
\special{pa 3062 3062}%
\special{fp}%
%
\special{pn 8}%
\special{pa 3072 950}%
\special{pa 3062 940}%
\special{fp}%
\special{sh 1}%
\special{pa 3062 940}%
\special{pa 3092 1002}%
\special{pa 3098 978}%
\special{pa 3122 976}%
\special{pa 3062 940}%
\special{fp}%
%
\special{pn 8}%
\special{sh 0}%
\special{ar 3500 2000 28 28  0.0000000 6.2831853}%
%
\special{pn 8}%
\special{sh 0.600}%
\special{ar 2000 500 28 28  0.0000000 6.2831853}%
%
\special{pn 8}%
\special{sh 0.600}%
\special{ar 2000 3500 28 28  0.0000000 6.2831853}%
%
\special{pn 8}%
\special{sh 0}%
\special{ar 500 2000 28 28  0.0000000 6.2831853}%
%
\special{pn 8}%
\special{pa 2000 1000}%
\special{pa 2000 630}%
\special{fp}%
\special{sh 1}%
\special{pa 2000 630}%
\special{pa 1980 698}%
\special{pa 2000 684}%
\special{pa 2020 698}%
\special{pa 2000 630}%
\special{fp}%
%
\special{pn 8}%
\special{pa 2000 3000}%
\special{pa 2000 3370}%
\special{fp}%
\special{sh 1}%
\special{pa 2000 3370}%
\special{pa 2020 3304}%
\special{pa 2000 3318}%
\special{pa 1980 3304}%
\special{pa 2000 3370}%
\special{fp}%
%
\special{pn 8}%
\special{pa 3360 2000}%
\special{pa 2990 2000}%
\special{fp}%
\special{sh 1}%
\special{pa 2990 2000}%
\special{pa 3058 2020}%
\special{pa 3044 2000}%
\special{pa 3058 1980}%
\special{pa 2990 2000}%
\special{fp}%
%
\special{pn 8}%
\special{pa 650 2000}%
\special{pa 1020 2000}%
\special{fp}%
\special{sh 1}%
\special{pa 1020 2000}%
\special{pa 954 1980}%
\special{pa 968 2000}%
\special{pa 954 2020}%
\special{pa 1020 2000}%
\special{fp}%
%
\special{pn 8}%
\special{pa 3400 1900}%
\special{pa 3138 1640}%
\special{fp}%
\special{sh 1}%
\special{pa 3138 1640}%
\special{pa 3172 1700}%
\special{pa 3176 1678}%
\special{pa 3200 1672}%
\special{pa 3138 1640}%
\special{fp}%
%
\special{pn 8}%
\special{pa 2350 850}%
\special{pa 2088 590}%
\special{fp}%
\special{sh 1}%
\special{pa 2088 590}%
\special{pa 2122 650}%
\special{pa 2126 628}%
\special{pa 2150 622}%
\special{pa 2088 590}%
\special{fp}%
%
\special{pn 8}%
\special{pa 610 1890}%
\special{pa 872 1630}%
\special{fp}%
\special{sh 1}%
\special{pa 872 1630}%
\special{pa 812 1662}%
\special{pa 834 1668}%
\special{pa 840 1690}%
\special{pa 872 1630}%
\special{fp}%
%
\special{pn 8}%
\special{pa 1640 860}%
\special{pa 1902 600}%
\special{fp}%
\special{sh 1}%
\special{pa 1902 600}%
\special{pa 1842 632}%
\special{pa 1864 638}%
\special{pa 1870 660}%
\special{pa 1902 600}%
\special{fp}%
%
\special{pn 8}%
\special{pa 2350 3150}%
\special{pa 2088 3412}%
\special{fp}%
\special{sh 1}%
\special{pa 2088 3412}%
\special{pa 2150 3378}%
\special{pa 2126 3374}%
\special{pa 2122 3350}%
\special{pa 2088 3412}%
\special{fp}%
%
\special{pn 8}%
\special{pa 3390 2110}%
\special{pa 3128 2372}%
\special{fp}%
\special{sh 1}%
\special{pa 3128 2372}%
\special{pa 3190 2338}%
\special{pa 3166 2334}%
\special{pa 3162 2310}%
\special{pa 3128 2372}%
\special{fp}%
%
\special{pn 8}%
\special{pa 1650 3150}%
\special{pa 1912 3412}%
\special{fp}%
\special{sh 1}%
\special{pa 1912 3412}%
\special{pa 1880 3350}%
\special{pa 1874 3374}%
\special{pa 1852 3378}%
\special{pa 1912 3412}%
\special{fp}%
%
\special{pn 8}%
\special{pa 610 2100}%
\special{pa 872 2362}%
\special{fp}%
\special{sh 1}%
\special{pa 872 2362}%
\special{pa 840 2300}%
\special{pa 834 2324}%
\special{pa 812 2328}%
\special{pa 872 2362}%
\special{fp}%
\put(35.0000,-23.0000){\makebox(0,0)[lt]{$\widehat{U}$}}%
\put(28.6000,-32.6000){\makebox(0,0)[lt]{$\widehat{f}$}}%
\end{picture}%
\caption{}
\label{fig8}
\end{center}
\end{figure}
One consequence of this is the famous Cartwright-Littlewood
fixed point theorem stated as Theorem \ref{21}
at the end of Section 4. Before giving the proof, 
we shall raise two examples of an invariant domain
with Carath\'eodory rotation number 0.

\begin{example} \label{example1}
There is a simple homeomorphism $h$ of $S^2$ which satisfies the
following conditions. 

(1) The homeomorphism $h$ preserves a continuum $X$.

(2) There is no periodic point in $X$

(3) $S^2\setminus X$ consists of three open discs $U_+$, $U_-$ and $V$.

(4) All the three open discs are invariant by $h$.

(5) The Carath\'eodory rotation number of $V$ is 0.

To construct $h$, we start with a Morse
Smale diffeomorphism
$g$ of the interval $[0,1]$ whose fixed points are $0$ and $1$.
Consider the suspension flow of $g$ on the annulus $S^1\times[0,1]$.
Define $h$ to be the time $\alpha$ map of the flow, where $\alpha$
is any irrational number. Choose one orbit $Y$ from $S^1\times(0,1)$
and let $X=S^1\times\{0,1\}\cup Y$ and $V=S^1\times[0,1]\setminus X$.
Finally extend $h$ to $S^2$ in an obvious way. See Figure \ref{fig9}.
\begin{figure}[htbp]
\begin{center}
\unitlength 0.1in
\begin{picture}( 30.0000, 30.0000)(  5.0000,-35.0000)
%
\special{pn 8}%
\special{ar 2000 2000 1500 1500  0.0000000 6.2831853}%
%
\special{pn 8}%
\special{ar 2000 2000 500 500  0.0000000 6.2831853}%
%
\special{pn 8}%
\special{pa 3430 2020}%
\special{pa 3428 2052}%
\special{pa 3424 2084}%
\special{pa 3420 2116}%
\special{pa 3414 2148}%
\special{pa 3408 2178}%
\special{pa 3402 2210}%
\special{pa 3392 2240}%
\special{pa 3384 2272}%
\special{pa 3374 2302}%
\special{pa 3364 2332}%
\special{pa 3354 2362}%
\special{pa 3344 2394}%
\special{pa 3334 2424}%
\special{pa 3324 2454}%
\special{pa 3314 2486}%
\special{pa 3304 2516}%
\special{pa 3292 2546}%
\special{pa 3282 2576}%
\special{pa 3268 2604}%
\special{pa 3256 2634}%
\special{pa 3240 2662}%
\special{pa 3226 2690}%
\special{pa 3210 2718}%
\special{pa 3192 2744}%
\special{pa 3176 2772}%
\special{pa 3158 2798}%
\special{pa 3140 2826}%
\special{pa 3130 2840}%
\special{sp 0.070}%
%
\special{pn 8}%
\special{pa 1520 1570}%
\special{pa 1500 1596}%
\special{pa 1478 1620}%
\special{pa 1458 1644}%
\special{pa 1438 1670}%
\special{pa 1418 1696}%
\special{pa 1400 1722}%
\special{pa 1384 1748}%
\special{pa 1368 1776}%
\special{pa 1354 1804}%
\special{pa 1342 1834}%
\special{pa 1332 1864}%
\special{pa 1322 1894}%
\special{pa 1316 1926}%
\special{pa 1310 1958}%
\special{pa 1306 1990}%
\special{pa 1304 2022}%
\special{pa 1302 2056}%
\special{pa 1302 2088}%
\special{pa 1304 2122}%
\special{pa 1306 2154}%
\special{pa 1310 2188}%
\special{pa 1316 2220}%
\special{pa 1322 2252}%
\special{pa 1328 2284}%
\special{pa 1338 2316}%
\special{pa 1348 2346}%
\special{pa 1358 2378}%
\special{pa 1372 2406}%
\special{pa 1386 2436}%
\special{pa 1402 2464}%
\special{pa 1418 2490}%
\special{pa 1436 2516}%
\special{pa 1458 2540}%
\special{pa 1478 2564}%
\special{pa 1502 2586}%
\special{pa 1526 2608}%
\special{pa 1552 2628}%
\special{pa 1578 2646}%
\special{pa 1604 2664}%
\special{pa 1632 2680}%
\special{pa 1662 2696}%
\special{pa 1690 2712}%
\special{pa 1720 2724}%
\special{pa 1750 2738}%
\special{pa 1780 2750}%
\special{pa 1810 2760}%
\special{pa 1842 2768}%
\special{pa 1872 2778}%
\special{pa 1904 2784}%
\special{pa 1936 2790}%
\special{pa 1968 2794}%
\special{pa 2000 2798}%
\special{pa 2032 2800}%
\special{pa 2064 2800}%
\special{pa 2096 2800}%
\special{pa 2128 2798}%
\special{pa 2160 2794}%
\special{pa 2192 2790}%
\special{pa 2224 2784}%
\special{pa 2256 2776}%
\special{pa 2288 2768}%
\special{pa 2318 2758}%
\special{pa 2350 2748}%
\special{pa 2380 2734}%
\special{pa 2410 2722}%
\special{pa 2438 2706}%
\special{pa 2468 2692}%
\special{pa 2496 2674}%
\special{pa 2522 2656}%
\special{pa 2548 2638}%
\special{pa 2574 2618}%
\special{pa 2598 2596}%
\special{pa 2622 2574}%
\special{pa 2646 2550}%
\special{pa 2666 2526}%
\special{pa 2688 2502}%
\special{pa 2708 2476}%
\special{pa 2726 2450}%
\special{pa 2744 2422}%
\special{pa 2760 2394}%
\special{pa 2776 2366}%
\special{pa 2790 2336}%
\special{pa 2802 2306}%
\special{pa 2814 2276}%
\special{pa 2824 2246}%
\special{pa 2832 2214}%
\special{pa 2840 2182}%
\special{pa 2848 2150}%
\special{pa 2852 2118}%
\special{pa 2856 2086}%
\special{pa 2860 2054}%
\special{pa 2860 2022}%
\special{pa 2860 1988}%
\special{pa 2860 1956}%
\special{pa 2858 1924}%
\special{pa 2854 1890}%
\special{pa 2848 1858}%
\special{pa 2842 1826}%
\special{pa 2834 1796}%
\special{pa 2826 1764}%
\special{pa 2816 1732}%
\special{pa 2804 1702}%
\special{pa 2792 1672}%
\special{pa 2778 1644}%
\special{pa 2762 1614}%
\special{pa 2746 1586}%
\special{pa 2730 1558}%
\special{pa 2710 1532}%
\special{pa 2692 1504}%
\special{pa 2672 1480}%
\special{pa 2650 1454}%
\special{pa 2628 1430}%
\special{pa 2606 1406}%
\special{pa 2582 1384}%
\special{pa 2558 1364}%
\special{pa 2534 1342}%
\special{pa 2508 1324}%
\special{pa 2482 1304}%
\special{pa 2454 1286}%
\special{pa 2426 1270}%
\special{pa 2398 1254}%
\special{pa 2370 1240}%
\special{pa 2342 1226}%
\special{pa 2312 1212}%
\special{pa 2282 1200}%
\special{pa 2252 1188}%
\special{pa 2222 1178}%
\special{pa 2190 1168}%
\special{pa 2160 1160}%
\special{pa 2128 1152}%
\special{pa 2096 1144}%
\special{pa 2064 1138}%
\special{pa 2032 1132}%
\special{pa 2000 1126}%
\special{pa 1968 1122}%
\special{pa 1936 1118}%
\special{pa 1904 1116}%
\special{pa 1872 1114}%
\special{pa 1840 1114}%
\special{pa 1808 1114}%
\special{pa 1776 1114}%
\special{pa 1744 1116}%
\special{pa 1712 1118}%
\special{pa 1680 1122}%
\special{pa 1648 1128}%
\special{pa 1616 1132}%
\special{pa 1584 1140}%
\special{pa 1552 1148}%
\special{pa 1522 1156}%
\special{pa 1490 1166}%
\special{pa 1460 1178}%
\special{pa 1430 1188}%
\special{pa 1400 1202}%
\special{pa 1370 1216}%
\special{pa 1342 1230}%
\special{pa 1314 1246}%
\special{pa 1286 1262}%
\special{pa 1258 1280}%
\special{pa 1230 1298}%
\special{pa 1204 1318}%
\special{pa 1180 1338}%
\special{pa 1154 1358}%
\special{pa 1130 1380}%
\special{pa 1108 1402}%
\special{pa 1086 1426}%
\special{pa 1064 1450}%
\special{pa 1044 1476}%
\special{pa 1024 1500}%
\special{pa 1004 1526}%
\special{pa 986 1554}%
\special{pa 970 1580}%
\special{pa 954 1608}%
\special{pa 938 1638}%
\special{pa 924 1666}%
\special{pa 910 1696}%
\special{pa 898 1726}%
\special{pa 886 1756}%
\special{pa 874 1786}%
\special{pa 864 1818}%
\special{pa 856 1850}%
\special{pa 846 1880}%
\special{pa 838 1912}%
\special{pa 832 1944}%
\special{pa 826 1976}%
\special{pa 820 2010}%
\special{pa 816 2042}%
\special{pa 812 2074}%
\special{pa 810 2108}%
\special{pa 808 2140}%
\special{pa 806 2174}%
\special{pa 806 2206}%
\special{pa 808 2238}%
\special{pa 808 2272}%
\special{pa 812 2304}%
\special{pa 814 2336}%
\special{pa 818 2368}%
\special{pa 824 2400}%
\special{pa 830 2432}%
\special{pa 836 2464}%
\special{pa 844 2496}%
\special{pa 852 2528}%
\special{pa 862 2558}%
\special{pa 872 2588}%
\special{pa 884 2618}%
\special{pa 896 2648}%
\special{pa 910 2678}%
\special{pa 924 2706}%
\special{pa 938 2736}%
\special{pa 954 2762}%
\special{pa 972 2790}%
\special{pa 990 2818}%
\special{pa 1008 2844}%
\special{pa 1028 2870}%
\special{pa 1048 2894}%
\special{pa 1068 2920}%
\special{pa 1090 2944}%
\special{pa 1112 2966}%
\special{pa 1136 2990}%
\special{pa 1160 3012}%
\special{pa 1184 3034}%
\special{pa 1208 3054}%
\special{pa 1234 3076}%
\special{pa 1260 3094}%
\special{pa 1286 3114}%
\special{pa 1314 3132}%
\special{pa 1342 3150}%
\special{pa 1370 3166}%
\special{pa 1398 3182}%
\special{pa 1426 3198}%
\special{pa 1456 3214}%
\special{pa 1484 3226}%
\special{pa 1514 3240}%
\special{pa 1544 3252}%
\special{pa 1574 3264}%
\special{pa 1604 3276}%
\special{pa 1636 3286}%
\special{pa 1666 3296}%
\special{pa 1698 3304}%
\special{pa 1728 3312}%
\special{pa 1760 3320}%
\special{pa 1792 3326}%
\special{pa 1822 3332}%
\special{pa 1854 3336}%
\special{pa 1886 3340}%
\special{pa 1918 3344}%
\special{pa 1950 3348}%
\special{pa 1982 3350}%
\special{pa 2014 3352}%
\special{pa 2046 3352}%
\special{pa 2078 3352}%
\special{pa 2110 3352}%
\special{pa 2142 3350}%
\special{pa 2174 3348}%
\special{pa 2206 3344}%
\special{pa 2238 3340}%
\special{pa 2270 3336}%
\special{pa 2302 3332}%
\special{pa 2334 3326}%
\special{pa 2366 3320}%
\special{pa 2398 3312}%
\special{pa 2428 3304}%
\special{pa 2460 3296}%
\special{pa 2490 3286}%
\special{pa 2522 3278}%
\special{pa 2552 3266}%
\special{pa 2582 3256}%
\special{pa 2612 3244}%
\special{pa 2642 3230}%
\special{pa 2672 3218}%
\special{pa 2702 3204}%
\special{pa 2730 3188}%
\special{pa 2758 3172}%
\special{pa 2786 3156}%
\special{pa 2814 3140}%
\special{pa 2842 3122}%
\special{pa 2868 3104}%
\special{pa 2894 3084}%
\special{pa 2918 3064}%
\special{pa 2944 3044}%
\special{pa 2966 3022}%
\special{pa 2990 3000}%
\special{pa 3012 2978}%
\special{pa 3034 2954}%
\special{pa 3056 2930}%
\special{pa 3076 2906}%
\special{pa 3096 2882}%
\special{pa 3116 2856}%
\special{pa 3136 2830}%
\special{pa 3154 2806}%
\special{pa 3174 2780}%
\special{pa 3180 2770}%
\special{sp}%
%
\special{pn 8}%
\special{pa 1510 1580}%
\special{pa 1534 1560}%
\special{pa 1560 1538}%
\special{pa 1586 1520}%
\special{pa 1612 1504}%
\special{pa 1642 1490}%
\special{pa 1672 1480}%
\special{pa 1702 1470}%
\special{pa 1734 1462}%
\special{pa 1766 1456}%
\special{pa 1798 1450}%
\special{pa 1828 1446}%
\special{pa 1860 1442}%
\special{pa 1892 1440}%
\special{pa 1924 1440}%
\special{pa 1956 1440}%
\special{pa 1988 1440}%
\special{pa 2000 1440}%
\special{sp 0.070}%
\put(31.3000,-19.3000){\makebox(0,0){$V$}}%
\end{picture}%
\caption{}
\label{fig9}
\end{center}
\end{figure}
Then the homeomorphism $\hat f$ on
the Carath\'eodory compactification $\hat V$ has two fixed prime ends.

\end{example}

\begin{example} \label{example2}
Let $g$ be a Denjoy's $C^1$ diffeomorphism of $S^1$ whose minimal set is a Cantor set $\mathfrak{N}$.
We put the suspension $T^2=S^1\times \R/(x,y)\sim (g(x),y+1)$.
For an irrational number $\alpha$, we define $f:T^2 \to T^2$ by $f([x,y])=[x+\alpha, y]$.
Then the minimal set of $f$ is $\mathfrak{N} \times \R/\sim$.
Its complement $U$ is a simply connected invariant domain.
By the same reason as in Example \ref{example1}, the
Carath\'eodory rotation number of $U$ is 0. See figure \ref{fig10}.
\begin{figure}[htbp]
\begin{center}
\unitlength 0.1in
%
\caption{}
\label{fig10}
\end{center}
\end{figure}
\end{example}

\bigskip
{\bf Proof of Theorem \ref{5}}. 
By the assumption on the Carath\'eodory rotation number,
the homeomorphism $\hat f$ has a fixed point $\xi$ in $\PP(U)$.
Let $\{c_i\}$ be an extendable topological chain 
representing $\xi$. Recall that
$\overline c_i$ are mutually disjoint in $\Sigma$. Also a ray
that is a half-ray in $c_i$ is
extendable and therefore belongs to some prime end by Lemma \ref{3}.
This implies that the cross cut $c_i$ is extendable in the
Carath\'eodory compactification $\hat U$. The closure of $c_i$
in $\hat U$ is denoted by $\hat c_i$.
By Lemma \ref{4a}
$\hat c_i$ are also mutually disjoint.

Assume for contradiction that $\hat f\hat c_i\cap \hat c_i\neq\emptyset$
for infinitely many $i$. Then again by Lemma \ref{4a} we have
$f\overline c_i\cap\overline c_i\neq\emptyset$.
Since diam$(c_i)\to0$, the point of accumulation of $c_i$ must be a fixed
point of $f$. 
Therefore we can assume $\hat f\hat c_i\cap \hat c_i=\emptyset$
for any $i$.

Let $\hat U(c_i)$ be the component of $\hat U\setminus \hat c_i$
not containing the base point $0\in U$. Notice that
$U(c_i)=U\cap\hat U(c_i)$.
Then we have for each large $i$ either $\hat f\hat c_i\subset \hat U(c_i)$ or
$\hat c_i\subset\hat f\hat U(c_i)$ because $\xi$ is a fixed point of $\hat f$.
Assume, to fix the idea, that
$\hat f\hat c_i\subset \hat U(c_i)$ for any $i$, by passing to a subsequence.

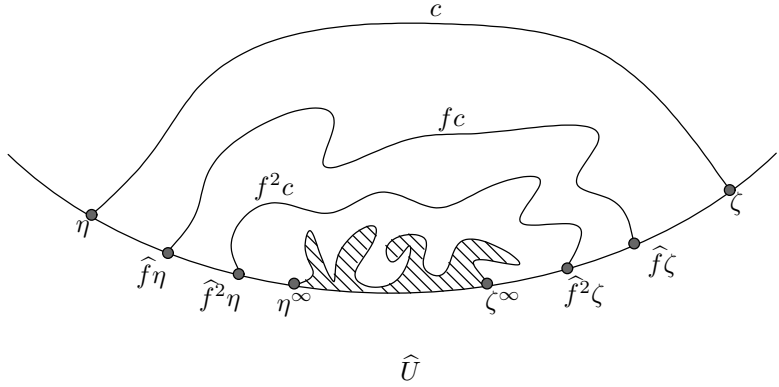
\begin{figure}[htbp]
\begin{center}
\unitlength 0.1in
\begin{picture}( 47.2300, 19.2500)( -2.1000,-23.0500)
%
\special{pn 8}%
\special{ar 2500 200 2500 1800  0.6344807 2.5025866}%
%
\special{pn 8}%
\special{pa 1990 1960}%
\special{pa 2028 1946}%
\special{pa 2062 1932}%
\special{pa 2084 1914}%
\special{pa 2092 1892}%
\special{pa 2084 1868}%
\special{pa 2068 1838}%
\special{pa 2054 1802}%
\special{pa 2046 1762}%
\special{pa 2050 1720}%
\special{pa 2058 1688}%
\special{pa 2070 1670}%
\special{pa 2082 1678}%
\special{pa 2094 1702}%
\special{pa 2108 1742}%
\special{pa 2122 1786}%
\special{pa 2140 1830}%
\special{pa 2156 1872}%
\special{pa 2172 1906}%
\special{pa 2188 1930}%
\special{pa 2198 1942}%
\special{pa 2202 1934}%
\special{pa 2204 1912}%
\special{pa 2206 1878}%
\special{pa 2210 1838}%
\special{pa 2222 1796}%
\special{pa 2244 1756}%
\special{pa 2274 1720}%
\special{pa 2306 1688}%
\special{pa 2340 1666}%
\special{pa 2372 1650}%
\special{pa 2398 1646}%
\special{pa 2416 1652}%
\special{pa 2422 1672}%
\special{pa 2416 1702}%
\special{pa 2402 1738}%
\special{pa 2382 1772}%
\special{pa 2360 1804}%
\special{pa 2340 1834}%
\special{pa 2322 1862}%
\special{pa 2310 1888}%
\special{pa 2308 1912}%
\special{pa 2318 1934}%
\special{pa 2340 1950}%
\special{pa 2370 1964}%
\special{pa 2406 1970}%
\special{pa 2446 1968}%
\special{pa 2486 1958}%
\special{pa 2522 1938}%
\special{pa 2552 1906}%
\special{pa 2576 1864}%
\special{pa 2592 1820}%
\special{pa 2600 1782}%
\special{pa 2600 1756}%
\special{pa 2592 1750}%
\special{pa 2576 1768}%
\special{pa 2552 1800}%
\special{pa 2528 1828}%
\special{pa 2502 1840}%
\special{pa 2482 1824}%
\special{pa 2470 1788}%
\special{pa 2476 1752}%
\special{pa 2500 1724}%
\special{pa 2536 1704}%
\special{pa 2576 1694}%
\special{pa 2618 1692}%
\special{pa 2654 1700}%
\special{pa 2678 1716}%
\special{pa 2684 1742}%
\special{pa 2680 1772}%
\special{pa 2670 1808}%
\special{pa 2664 1842}%
\special{pa 2672 1874}%
\special{pa 2696 1896}%
\special{pa 2728 1910}%
\special{pa 2758 1908}%
\special{pa 2780 1892}%
\special{pa 2796 1866}%
\special{pa 2808 1832}%
\special{pa 2820 1796}%
\special{pa 2834 1764}%
\special{pa 2854 1736}%
\special{pa 2880 1720}%
\special{pa 2914 1716}%
\special{pa 2954 1720}%
\special{pa 2996 1732}%
\special{pa 3038 1750}%
\special{pa 3080 1770}%
\special{pa 3116 1790}%
\special{pa 3144 1812}%
\special{pa 3164 1830}%
\special{pa 3172 1844}%
\special{pa 3168 1852}%
\special{pa 3148 1852}%
\special{pa 3118 1848}%
\special{pa 3080 1842}%
\special{pa 3038 1834}%
\special{pa 2996 1826}%
\special{pa 2960 1824}%
\special{pa 2934 1828}%
\special{pa 2926 1842}%
\special{pa 2936 1864}%
\special{pa 2958 1892}%
\special{pa 2986 1926}%
\special{pa 3010 1950}%
\special{sp}%
%
\special{pn 8}%
\special{pa 1710 1900}%
\special{pa 1696 1870}%
\special{pa 1682 1840}%
\special{pa 1670 1810}%
\special{pa 1662 1778}%
\special{pa 1660 1748}%
\special{pa 1664 1716}%
\special{pa 1674 1686}%
\special{pa 1688 1654}%
\special{pa 1708 1626}%
\special{pa 1730 1600}%
\special{pa 1756 1576}%
\special{pa 1786 1556}%
\special{pa 1816 1540}%
\special{pa 1848 1532}%
\special{pa 1880 1526}%
\special{pa 1912 1528}%
\special{pa 1944 1532}%
\special{pa 1976 1540}%
\special{pa 2008 1550}%
\special{pa 2040 1560}%
\special{pa 2072 1568}%
\special{pa 2102 1576}%
\special{pa 2134 1582}%
\special{pa 2164 1582}%
\special{pa 2196 1578}%
\special{pa 2224 1568}%
\special{pa 2254 1556}%
\special{pa 2284 1540}%
\special{pa 2312 1524}%
\special{pa 2342 1508}%
\special{pa 2372 1492}%
\special{pa 2402 1480}%
\special{pa 2432 1472}%
\special{pa 2462 1470}%
\special{pa 2492 1474}%
\special{pa 2524 1482}%
\special{pa 2556 1492}%
\special{pa 2586 1506}%
\special{pa 2618 1518}%
\special{pa 2648 1532}%
\special{pa 2678 1542}%
\special{pa 2706 1550}%
\special{pa 2734 1554}%
\special{pa 2762 1552}%
\special{pa 2788 1548}%
\special{pa 2814 1538}%
\special{pa 2842 1528}%
\special{pa 2868 1516}%
\special{pa 2896 1500}%
\special{pa 2926 1486}%
\special{pa 2958 1468}%
\special{pa 2992 1452}%
\special{pa 3028 1436}%
\special{pa 3066 1422}%
\special{pa 3108 1408}%
\special{pa 3152 1398}%
\special{pa 3192 1390}%
\special{pa 3228 1388}%
\special{pa 3256 1390}%
\special{pa 3276 1400}%
\special{pa 3280 1418}%
\special{pa 3274 1442}%
\special{pa 3256 1470}%
\special{pa 3234 1504}%
\special{pa 3210 1536}%
\special{pa 3190 1568}%
\special{pa 3174 1596}%
\special{pa 3170 1620}%
\special{pa 3180 1634}%
\special{pa 3200 1642}%
\special{pa 3232 1644}%
\special{pa 3270 1644}%
\special{pa 3310 1640}%
\special{pa 3354 1638}%
\special{pa 3394 1636}%
\special{pa 3432 1636}%
\special{pa 3462 1642}%
\special{pa 3484 1654}%
\special{pa 3494 1672}%
\special{pa 3494 1696}%
\special{pa 3486 1726}%
\special{pa 3474 1760}%
\special{pa 3456 1798}%
\special{pa 3434 1836}%
\special{pa 3420 1860}%
\special{sp}%
%
\special{pn 8}%
\special{pa 1330 1780}%
\special{pa 1352 1756}%
\special{pa 1374 1730}%
\special{pa 1396 1706}%
\special{pa 1416 1680}%
\special{pa 1436 1654}%
\special{pa 1452 1628}%
\special{pa 1468 1600}%
\special{pa 1480 1572}%
\special{pa 1490 1544}%
\special{pa 1498 1514}%
\special{pa 1504 1484}%
\special{pa 1512 1454}%
\special{pa 1520 1422}%
\special{pa 1530 1392}%
\special{pa 1544 1360}%
\special{pa 1562 1330}%
\special{pa 1584 1298}%
\special{pa 1610 1268}%
\special{pa 1640 1240}%
\special{pa 1674 1210}%
\special{pa 1712 1184}%
\special{pa 1750 1158}%
\special{pa 1790 1134}%
\special{pa 1832 1112}%
\special{pa 1874 1092}%
\special{pa 1916 1074}%
\special{pa 1956 1058}%
\special{pa 1996 1046}%
\special{pa 2034 1038}%
\special{pa 2070 1032}%
\special{pa 2102 1030}%
\special{pa 2132 1032}%
\special{pa 2156 1038}%
\special{pa 2176 1048}%
\special{pa 2190 1064}%
\special{pa 2198 1082}%
\special{pa 2200 1108}%
\special{pa 2196 1136}%
\special{pa 2186 1170}%
\special{pa 2174 1204}%
\special{pa 2160 1240}%
\special{pa 2150 1272}%
\special{pa 2142 1300}%
\special{pa 2140 1322}%
\special{pa 2148 1336}%
\special{pa 2162 1342}%
\special{pa 2182 1344}%
\special{pa 2210 1340}%
\special{pa 2240 1330}%
\special{pa 2276 1318}%
\special{pa 2314 1304}%
\special{pa 2352 1286}%
\special{pa 2394 1268}%
\special{pa 2434 1252}%
\special{pa 2474 1234}%
\special{pa 2514 1220}%
\special{pa 2550 1206}%
\special{pa 2584 1196}%
\special{pa 2616 1188}%
\special{pa 2646 1182}%
\special{pa 2676 1176}%
\special{pa 2704 1174}%
\special{pa 2732 1170}%
\special{pa 2762 1170}%
\special{pa 2790 1168}%
\special{pa 2820 1168}%
\special{pa 2850 1166}%
\special{pa 2882 1166}%
\special{pa 2916 1166}%
\special{pa 2950 1164}%
\special{pa 2988 1160}%
\special{pa 3030 1156}%
\special{pa 3072 1152}%
\special{pa 3118 1148}%
\special{pa 3164 1142}%
\special{pa 3210 1136}%
\special{pa 3256 1132}%
\special{pa 3302 1126}%
\special{pa 3346 1124}%
\special{pa 3388 1120}%
\special{pa 3428 1120}%
\special{pa 3464 1120}%
\special{pa 3498 1124}%
\special{pa 3526 1128}%
\special{pa 3550 1136}%
\special{pa 3568 1146}%
\special{pa 3580 1160}%
\special{pa 3586 1176}%
\special{pa 3584 1198}%
\special{pa 3576 1222}%
\special{pa 3558 1250}%
\special{pa 3538 1282}%
\special{pa 3516 1314}%
\special{pa 3496 1344}%
\special{pa 3478 1374}%
\special{pa 3468 1402}%
\special{pa 3468 1426}%
\special{pa 3480 1444}%
\special{pa 3504 1458}%
\special{pa 3534 1470}%
\special{pa 3568 1480}%
\special{pa 3604 1492}%
\special{pa 3638 1504}%
\special{pa 3668 1522}%
\special{pa 3692 1542}%
\special{pa 3712 1566}%
\special{pa 3728 1594}%
\special{pa 3742 1624}%
\special{pa 3750 1656}%
\special{pa 3758 1690}%
\special{pa 3766 1724}%
\special{pa 3770 1750}%
\special{sp}%
%
\special{pn 8}%
\special{pa 930 1590}%
\special{pa 954 1568}%
\special{pa 978 1548}%
\special{pa 1002 1526}%
\special{pa 1024 1504}%
\special{pa 1048 1482}%
\special{pa 1070 1460}%
\special{pa 1094 1436}%
\special{pa 1116 1414}%
\special{pa 1138 1390}%
\special{pa 1160 1368}%
\special{pa 1182 1344}%
\special{pa 1202 1318}%
\special{pa 1222 1294}%
\special{pa 1242 1268}%
\special{pa 1262 1242}%
\special{pa 1280 1216}%
\special{pa 1298 1190}%
\special{pa 1316 1162}%
\special{pa 1334 1136}%
\special{pa 1352 1108}%
\special{pa 1370 1082}%
\special{pa 1388 1054}%
\special{pa 1406 1028}%
\special{pa 1424 1002}%
\special{pa 1444 978}%
\special{pa 1462 952}%
\special{pa 1482 928}%
\special{pa 1504 906}%
\special{pa 1526 884}%
\special{pa 1548 864}%
\special{pa 1572 844}%
\special{pa 1596 824}%
\special{pa 1620 806}%
\special{pa 1646 790}%
\special{pa 1670 774}%
\special{pa 1698 758}%
\special{pa 1724 744}%
\special{pa 1752 730}%
\special{pa 1780 716}%
\special{pa 1810 704}%
\special{pa 1838 692}%
\special{pa 1868 682}%
\special{pa 1898 672}%
\special{pa 1930 662}%
\special{pa 1960 654}%
\special{pa 1992 646}%
\special{pa 2024 638}%
\special{pa 2056 632}%
\special{pa 2090 626}%
\special{pa 2122 620}%
\special{pa 2156 614}%
\special{pa 2190 610}%
\special{pa 2222 606}%
\special{pa 2256 602}%
\special{pa 2292 598}%
\special{pa 2326 596}%
\special{pa 2360 592}%
\special{pa 2394 590}%
\special{pa 2430 590}%
\special{pa 2464 588}%
\special{pa 2500 588}%
\special{pa 2534 586}%
\special{pa 2570 586}%
\special{pa 2604 586}%
\special{pa 2640 586}%
\special{pa 2674 586}%
\special{pa 2710 588}%
\special{pa 2744 588}%
\special{pa 2780 590}%
\special{pa 2814 590}%
\special{pa 2848 592}%
\special{pa 2882 594}%
\special{pa 2918 594}%
\special{pa 2950 596}%
\special{pa 2984 598}%
\special{pa 3018 602}%
\special{pa 3052 604}%
\special{pa 3084 606}%
\special{pa 3118 610}%
\special{pa 3150 614}%
\special{pa 3182 618}%
\special{pa 3214 624}%
\special{pa 3246 628}%
\special{pa 3278 634}%
\special{pa 3308 640}%
\special{pa 3338 648}%
\special{pa 3370 656}%
\special{pa 3398 664}%
\special{pa 3428 672}%
\special{pa 3458 682}%
\special{pa 3486 694}%
\special{pa 3514 704}%
\special{pa 3542 718}%
\special{pa 3570 730}%
\special{pa 3596 744}%
\special{pa 3624 758}%
\special{pa 3650 774}%
\special{pa 3674 792}%
\special{pa 3700 808}%
\special{pa 3724 826}%
\special{pa 3748 846}%
\special{pa 3772 864}%
\special{pa 3796 886}%
\special{pa 3820 906}%
\special{pa 3842 928}%
\special{pa 3866 950}%
\special{pa 3888 972}%
\special{pa 3910 996}%
\special{pa 3932 1018}%
\special{pa 3954 1044}%
\special{pa 3974 1068}%
\special{pa 3996 1092}%
\special{pa 4016 1118}%
\special{pa 4038 1144}%
\special{pa 4058 1170}%
\special{pa 4078 1196}%
\special{pa 4100 1222}%
\special{pa 4120 1250}%
\special{pa 4140 1278}%
\special{pa 4160 1304}%
\special{pa 4180 1332}%
\special{pa 4198 1360}%
\special{pa 4218 1388}%
\special{pa 4238 1414}%
\special{pa 4258 1442}%
\special{pa 4270 1460}%
\special{sp}%
\put(9.2000,-16.2000){\makebox(0,0)[rt]{$\eta$}}%
\put(13.2000,-18.1000){\makebox(0,0)[rt]{$\widehat{f}\eta$}}%
\put(17.1000,-19.3000){\makebox(0,0)[rt]{$\widehat{f}^2\eta$}}%
\put(20.9000,-20.0000){\makebox(0,0)[rt]{$\eta^\infty$}}%
\put(29.9000,-19.9000){\makebox(0,0)[lt]{$\zeta^\infty$}}%
\put(34.0000,-19.1000){\makebox(0,0)[lt]{$\widehat{f}^2\zeta$}}%
\put(38.4000,-17.5000){\makebox(0,0)[lt]{$\widehat{f}\zeta$}}%
\put(42.7000,-14.7000){\makebox(0,0)[lt]{$\zeta$}}%
\put(27.0000,-5.5000){\makebox(0,0)[lb]{$c$}}%
\put(27.4000,-11.5000){\makebox(0,0)[lb]{$fc$}}%
\put(17.8000,-15.1000){\makebox(0,0)[lb]{$f^2c$}}%
%
\special{pn 8}%
\special{sh 0.600}%
\special{ar 4270 1460 28 28  0.0000000 6.2831853}%
%
\special{pn 8}%
\special{sh 0.600}%
\special{ar 930 1590 28 28  0.0000000 6.2831853}%
%
\special{pn 8}%
\special{sh 0.600}%
\special{ar 1330 1790 28 28  0.0000000 6.2831853}%
%
\special{pn 8}%
\special{sh 0.600}%
\special{ar 1700 1900 28 28  0.0000000 6.2831853}%
%
\special{pn 8}%
\special{sh 0.600}%
\special{ar 1990 1950 28 28  0.0000000 6.2831853}%
%
\special{pn 8}%
\special{sh 0.600}%
\special{ar 3000 1950 28 28  0.0000000 6.2831853}%
%
\special{pn 8}%
\special{sh 0.600}%
\special{ar 3420 1870 28 28  0.0000000 6.2831853}%
%
\special{pn 8}%
\special{sh 0.600}%
\special{ar 3770 1740 28 28  0.0000000 6.2831853}%
\put(26.0000,-23.9000){\makebox(0,0){$\widehat{U}$}}%
%
\special{pn 8}%
\special{pa 2930 1970}%
\special{pa 2800 1840}%
\special{fp}%
\special{pa 2970 1950}%
\special{pa 2820 1800}%
\special{fp}%
\special{pa 2930 1850}%
\special{pa 2840 1760}%
\special{fp}%
\special{pa 2960 1820}%
\special{pa 2870 1730}%
\special{fp}%
\special{pa 3030 1830}%
\special{pa 2920 1720}%
\special{fp}%
\special{pa 3100 1840}%
\special{pa 2990 1730}%
\special{fp}%
\special{pa 2880 1980}%
\special{pa 2790 1890}%
\special{fp}%
\special{pa 2820 1980}%
\special{pa 2750 1910}%
\special{fp}%
\special{pa 2770 1990}%
\special{pa 2590 1810}%
\special{fp}%
\special{pa 2710 1990}%
\special{pa 2580 1860}%
\special{fp}%
\special{pa 2650 1990}%
\special{pa 2560 1900}%
\special{fp}%
\special{pa 2600 2000}%
\special{pa 2530 1930}%
\special{fp}%
\special{pa 2540 2000}%
\special{pa 2490 1950}%
\special{fp}%
\special{pa 2480 2000}%
\special{pa 2450 1970}%
\special{fp}%
\special{pa 2660 1820}%
\special{pa 2600 1760}%
\special{fp}%
\special{pa 2670 1770}%
\special{pa 2590 1690}%
\special{fp}%
\special{pa 2590 1750}%
\special{pa 2540 1700}%
\special{fp}%
\special{pa 2560 1780}%
\special{pa 2500 1720}%
\special{fp}%
\special{pa 2530 1810}%
\special{pa 2470 1750}%
\special{fp}%
%
\special{pn 8}%
\special{pa 2310 1890}%
\special{pa 2220 1800}%
\special{fp}%
\special{pa 2320 1840}%
\special{pa 2240 1760}%
\special{fp}%
\special{pa 2350 1810}%
\special{pa 2270 1730}%
\special{fp}%
\special{pa 2370 1770}%
\special{pa 2300 1700}%
\special{fp}%
\special{pa 2400 1740}%
\special{pa 2330 1670}%
\special{fp}%
\special{pa 2420 1700}%
\special{pa 2370 1650}%
\special{fp}%
\special{pa 2350 1990}%
\special{pa 2210 1850}%
\special{fp}%
\special{pa 2290 1990}%
\special{pa 2200 1900}%
\special{fp}%
\special{pa 2230 1990}%
\special{pa 2060 1820}%
\special{fp}%
\special{pa 2140 1840}%
\special{pa 2050 1750}%
\special{fp}%
\special{pa 2110 1750}%
\special{pa 2060 1700}%
\special{fp}%
\special{pa 2160 1980}%
\special{pa 2090 1910}%
\special{fp}%
\special{pa 2090 1970}%
\special{pa 2050 1930}%
\special{fp}%
\end{picture}%
\caption{Hatched area is $U\setminus U_0$}
\label{fig11}
\end{center}
\end{figure}

Now let $N$ be an neighbourhood of the frontier ${\rm Fr}(U)$ 
which does not intersect the fixed point set ${\rm Fix}(f)$
of $f$. Then since 
$\cap_i\overline U(c_i)\subset {\rm Fr}(U)$ in $\Sigma$ by Lemma \ref{0a}, the
closure of the domain $U(c_i)$ for some big $i$ is contained in $N$.
Fix once and for all such $c_i$ and denote it by $c$.
The two end points $\eta$ and $\zeta$ of $\hat c$ form an interval
$[\eta,\zeta]$ in $\PP(U)$ containing the prime end $\xi$, a fixed point
of $\hat f$.
On this interval we have 
$$
\eta<\hat f\eta<\hat f^2\eta<\cdots<\hat f^2\zeta<\hat f\zeta<\zeta.
$$
Assume that 
\begin{equation} \label{6}
\eta^\infty=\lim \hat f^n\eta<\zeta^\infty=\lim\hat f^n\zeta.
\end{equation}
See Figure \ref{fig11}.
A contradiction will show that the map $\hat f$ is Morse-Smale
on $\PP(U)$.

Consider a domain 
$$U_0=U\setminus\cap_nf^nU(c)
$$
and notice that ${\rm Fix(f)}\cap{\rm Fr}U_0=\emptyset$,
by the choice of $c$. 
The chain $\{f^nc\}$ of $U$ is
also a chain of $U_0$, and each cross cut $f^nc$ is of course
extendable. 
An important feature of $U_0$ is that the intersection
of the contents is empty, i.\ e.\
\begin{equation} \label{7}
\cap_{n=0}^\infty f^nU_0(c)=\emptyset.
\end{equation}

Let us denote by $\hat f_0$ the homeomorphism
induced by $f$ on the Carath\'eodory compactification
$\hat U_0$  of $U_0$. Let $\eta_0$ and $\zeta_0$ be
the prime ends in $\PP(U_0)$ corresponding to the end points of
$c$. Then we have
$$
\eta_0<\hat f_0\eta_0<\hat f_0^2\eta_0<\cdots<\hat f_0^2\zeta_0<\hat f_0
\zeta_0<\zeta_0.
$$
Let 
$\eta_0^\infty=\lim \hat f_0^n\eta_0$ and
$\zeta_0^\infty=\lim\hat f_0^n\zeta_0$.
It follows from the definition of topological chains
that there is an order preserving homeomorphism between 
$\PP(U)\setminus[\eta^\infty,\zeta^\infty]$
and
$\PP(U_0)\setminus[\eta_0^\infty,\zeta_0^\infty]$.
Let us show $\eta_0^\infty<\zeta_0^\infty$.
Assuming the contrary, we get an extendable topological chain
$c'_i$ representing $\eta_0^\infty=\zeta_0^\infty$.
Let $\alpha_0^i$ and $\beta_0^i$ be the two prime ends in $\PP(U_0)$
corresponding to $c'_i$. Then clearly the sequences $\hat f_0^n\eta_0$
and $\alpha_0^i$ have the same limit $\eta_0^\infty=\zeta_0^\infty$. 
In other words, they are cofinal,
that is,
for any $i$, there is $n$ such that $\alpha_0^i<\hat f_0^n\eta_0$
and for any $n$, there is $i$ such that $\hat f_0^n\eta_0<\alpha_0^i$.
Likewise $\beta_0^i$ and $\hat f_0^n\zeta$ are cofinal.
Now $c'_i$ is also an extendable topological chain of $U$
joining $\alpha_i$ and $\beta_i$ in $\PP(U)$. 
Since $\PP(U)\setminus[\eta^\infty,\zeta^\infty]$
and
$\PP(U_0)\setminus[\eta_0^\infty,\zeta_0^\infty]$
are order preserving homeomorphic, we see that
$\alpha_i$ and $\hat f^n\eta$ are cofinal
and $\beta_i$ and $\hat f^n\zeta$ are cofinal.
Since $\{c'_i\}$ is also a topological chain of $U$,
this shows that $\eta^\infty=\zeta^\infty$,
against the assumption (\ref{6}).

Since $f$ is fixed point free on ${\rm Fr}(U_0)$
and the natural map $\Phi:\EE(U_0)\to {\rm Fr}(U_0)$
is equivariant, $\Phi\circ\hat f_0=f\circ\Phi$,
the set of extendable ends $\EE(U_0)$ 
is disjoint from  ${\rm Fix}(\hat f_0)$. Lemma \ref{4} implies that
the fixed point set of $\hat f_0$ is nowhere dense in $\PP(U_0)$.
Thus there is a point $\sigma$ in the interval 
$[\eta_0^\infty,\zeta_0^\infty]$ which is not fixed by $\hat f_0$.
See Figure \ref{fig13}.
\begin{figure}[htbp]
\begin{center}
\unitlength 0.1in
\begin{picture}( 50.0000, 10.9000)(  5.0000,-25.9000)
%
\special{pn 8}%
\special{pa 500 2500}%
\special{pa 5500 2500}%
\special{fp}%
%
\special{pn 8}%
\special{ar 3000 2500 2000 1000  3.1415927 6.2831853}%
%
\special{pn 8}%
\special{ar 3000 2500 1500 700  3.1415927 6.2831853}%
%
\special{pn 8}%
\special{ar 3000 2510 230 230  3.1891758 6.2397344}%
%
\special{pn 8}%
\special{ar 3000 2500 400 400  3.1415927 6.2831853}%
%
\special{pn 8}%
\special{sh 0.600}%
\special{ar 4200 2500 28 28  0.0000000 6.2831853}%
%
\special{pn 8}%
\special{sh 0.600}%
\special{ar 1810 2500 28 28  0.0000000 6.2831853}%
%
\special{pn 8}%
\special{sh 0.600}%
\special{ar 3000 2500 28 28  0.0000000 6.2831853}%
%
\special{pn 8}%
\special{pa 500 2500}%
\special{pa 800 2500}%
\special{fp}%
\special{sh 1}%
\special{pa 800 2500}%
\special{pa 734 2480}%
\special{pa 748 2500}%
\special{pa 734 2520}%
\special{pa 800 2500}%
\special{fp}%
%
\special{pn 8}%
\special{pa 3000 2500}%
\special{pa 3150 2500}%
\special{fp}%
\special{sh 1}%
\special{pa 3150 2500}%
\special{pa 3084 2480}%
\special{pa 3098 2500}%
\special{pa 3084 2520}%
\special{pa 3150 2500}%
\special{fp}%
%
\special{pn 8}%
\special{pa 3230 2500}%
\special{pa 3350 2500}%
\special{fp}%
\special{sh 1}%
\special{pa 3350 2500}%
\special{pa 3284 2480}%
\special{pa 3298 2500}%
\special{pa 3284 2520}%
\special{pa 3350 2500}%
\special{fp}%
%
\special{pn 8}%
\special{pa 3000 1540}%
\special{pa 3000 1780}%
\special{fp}%
\special{sh 1}%
\special{pa 3000 1780}%
\special{pa 3020 1714}%
\special{pa 3000 1728}%
\special{pa 2980 1714}%
\special{pa 3000 1780}%
\special{fp}%
%
\special{pn 8}%
\special{pa 4590 1940}%
\special{pa 4250 2060}%
\special{fp}%
\special{sh 1}%
\special{pa 4250 2060}%
\special{pa 4320 2058}%
\special{pa 4300 2042}%
\special{pa 4306 2020}%
\special{pa 4250 2060}%
\special{fp}%
%
\special{pn 8}%
\special{pa 1400 1940}%
\special{pa 1740 2060}%
\special{fp}%
\special{sh 1}%
\special{pa 1740 2060}%
\special{pa 1684 2020}%
\special{pa 1690 2042}%
\special{pa 1670 2058}%
\special{pa 1740 2060}%
\special{fp}%
%
\special{pn 8}%
\special{pa 3860 1640}%
\special{pa 3660 1840}%
\special{fp}%
\special{sh 1}%
\special{pa 3660 1840}%
\special{pa 3722 1808}%
\special{pa 3698 1802}%
\special{pa 3694 1780}%
\special{pa 3660 1840}%
\special{fp}%
%
\special{pn 8}%
\special{pa 2140 1640}%
\special{pa 2340 1840}%
\special{fp}%
\special{sh 1}%
\special{pa 2340 1840}%
\special{pa 2308 1780}%
\special{pa 2302 1802}%
\special{pa 2280 1808}%
\special{pa 2340 1840}%
\special{fp}%
\put(30.3000,-17.3000){\makebox(0,0)[lb]{$\widehat{f}$}}%
\put(17.9000,-25.6000){\makebox(0,0)[lt]{$\eta^\infty_0$}}%
%
\special{pn 8}%
\special{pa 3180 2330}%
\special{pa 3270 2240}%
\special{fp}%
\special{sh 1}%
\special{pa 3270 2240}%
\special{pa 3210 2274}%
\special{pa 3232 2278}%
\special{pa 3238 2302}%
\special{pa 3270 2240}%
\special{fp}%
%
\special{pn 8}%
\special{pa 2820 2320}%
\special{pa 2730 2230}%
\special{fp}%
\special{sh 1}%
\special{pa 2730 2230}%
\special{pa 2764 2292}%
\special{pa 2768 2268}%
\special{pa 2792 2264}%
\special{pa 2730 2230}%
\special{fp}%
\put(29.5000,-25.7000){\makebox(0,0)[lt]{$\tau$}}%
\put(35.5000,-25.8000){\makebox(0,0)[lt]{$\sigma$}}%
\put(41.8000,-25.9000){\makebox(0,0)[lt]{$\zeta^\infty_0$}}%
%
\special{pn 8}%
\special{sh 0.600}%
\special{ar 3600 2500 28 28  0.0000000 6.2831853}%
%
\special{pn 8}%
\special{pa 3400 2500}%
\special{pa 3500 2500}%
\special{fp}%
\special{sh 1}%
\special{pa 3500 2500}%
\special{pa 3434 2480}%
\special{pa 3448 2500}%
\special{pa 3434 2520}%
\special{pa 3500 2500}%
\special{fp}%
%
\special{pn 8}%
\special{pa 3000 2270}%
\special{pa 3000 2120}%
\special{fp}%
\special{sh 1}%
\special{pa 3000 2120}%
\special{pa 2980 2188}%
\special{pa 3000 2174}%
\special{pa 3020 2188}%
\special{pa 3000 2120}%
\special{fp}%
%
\special{pn 8}%
\special{pa 1000 2500}%
\special{pa 1300 2500}%
\special{fp}%
\special{sh 1}%
\special{pa 1300 2500}%
\special{pa 1234 2480}%
\special{pa 1248 2500}%
\special{pa 1234 2520}%
\special{pa 1300 2500}%
\special{fp}%
%
\special{pn 8}%
\special{pa 5000 2500}%
\special{pa 4750 2500}%
\special{fp}%
\special{sh 1}%
\special{pa 4750 2500}%
\special{pa 4818 2520}%
\special{pa 4804 2500}%
\special{pa 4818 2480}%
\special{pa 4750 2500}%
\special{fp}%
%
\special{pn 8}%
\special{pa 5500 2500}%
\special{pa 5250 2500}%
\special{fp}%
\special{sh 1}%
\special{pa 5250 2500}%
\special{pa 5318 2520}%
\special{pa 5304 2500}%
\special{pa 5318 2480}%
\special{pa 5250 2500}%
\special{fp}%
%
\special{pn 8}%
\special{pa 3700 2500}%
\special{pa 3800 2500}%
\special{fp}%
\special{sh 1}%
\special{pa 3800 2500}%
\special{pa 3734 2480}%
\special{pa 3748 2500}%
\special{pa 3734 2520}%
\special{pa 3800 2500}%
\special{fp}%
\end{picture}%
\caption{}
\label{fig13}
\end{center}
\end{figure}
To fix the idea assume $\hat f_0\sigma>\sigma$ and let
$\hat f_0^{-n}\sigma\downarrow \tau$. Let $\{c''_i\}$ be an extendable
topological chain of $U_0$ representing $\tau$. Denote
by $U_0(c''_i)$ the content of $c''_i$ in $U_0$.
As before we have $\hat f_0(c''_i)\cap c''_i=\emptyset$ if we pass
to a subsequence. But $\tau$ is repelling on its right side.
Therefore $U_0(c''_i)\subset\hat f_0U_0(c''_i)$.
If we choose $i$ big enough, we have $U_0(c''_i)\subset U_0(c)$.
But this is contrary to (\ref{7}), finishing the proof
that $\hat f$ is Morse-Smale on $\PP(U)$.

\begin{figure}[htbp]
\begin{center}
\unitlength 0.1in
\begin{picture}( 59.9400, 16.2000)(  8.0000,-20.0000)
%
\special{pn 8}%
\special{pa 2000 2000}%
\special{pa 800 800}%
\special{fp}%
%
\special{pn 8}%
\special{pa 2000 2000}%
\special{pa 3200 800}%
\special{fp}%
%
\special{pn 20}%
\special{pa 2000 2000}%
\special{pa 2000 1968}%
\special{pa 1998 1934}%
\special{pa 1998 1902}%
\special{pa 1996 1868}%
\special{pa 1996 1836}%
\special{pa 1996 1802}%
\special{pa 1998 1770}%
\special{pa 2000 1738}%
\special{pa 2004 1706}%
\special{pa 2008 1676}%
\special{pa 2014 1644}%
\special{pa 2020 1614}%
\special{pa 2030 1584}%
\special{pa 2040 1554}%
\special{pa 2052 1524}%
\special{pa 2064 1496}%
\special{pa 2080 1468}%
\special{pa 2094 1440}%
\special{pa 2112 1412}%
\special{pa 2128 1384}%
\special{pa 2146 1358}%
\special{pa 2150 1350}%
\special{sp}%
%
\special{pn 20}%
\special{pa 2000 2000}%
\special{pa 2002 1970}%
\special{pa 2002 1938}%
\special{pa 2002 1906}%
\special{pa 2002 1874}%
\special{pa 2002 1842}%
\special{pa 2002 1810}%
\special{pa 2000 1776}%
\special{pa 1998 1742}%
\special{pa 1992 1708}%
\special{pa 1986 1674}%
\special{pa 1976 1644}%
\special{pa 1962 1618}%
\special{pa 1944 1594}%
\special{pa 1922 1576}%
\special{pa 1896 1560}%
\special{pa 1866 1548}%
\special{pa 1834 1538}%
\special{pa 1800 1530}%
\special{pa 1766 1522}%
\special{pa 1760 1520}%
\special{sp}%
%
\special{pn 20}%
\special{pa 2000 2000}%
\special{pa 2000 1970}%
\special{pa 2000 1938}%
\special{pa 2002 1906}%
\special{pa 2004 1876}%
\special{pa 2008 1842}%
\special{pa 2014 1810}%
\special{pa 2022 1778}%
\special{pa 2032 1744}%
\special{pa 2046 1712}%
\special{pa 2062 1682}%
\special{pa 2084 1660}%
\special{pa 2108 1642}%
\special{pa 2136 1632}%
\special{pa 2168 1630}%
\special{pa 2202 1630}%
\special{pa 2210 1630}%
\special{sp}%
%
\special{pn 20}%
\special{pa 2000 2000}%
\special{pa 2000 1968}%
\special{pa 2002 1936}%
\special{pa 2002 1904}%
\special{pa 2002 1872}%
\special{pa 2002 1840}%
\special{pa 2000 1808}%
\special{pa 2000 1776}%
\special{pa 2000 1744}%
\special{pa 1998 1712}%
\special{pa 1996 1680}%
\special{pa 1994 1648}%
\special{pa 1992 1616}%
\special{pa 1988 1584}%
\special{pa 1982 1552}%
\special{pa 1976 1522}%
\special{pa 1968 1490}%
\special{pa 1958 1460}%
\special{pa 1948 1430}%
\special{pa 1938 1400}%
\special{pa 1926 1370}%
\special{pa 1912 1340}%
\special{pa 1900 1312}%
\special{pa 1886 1282}%
\special{pa 1872 1254}%
\special{pa 1870 1250}%
\special{sp}%
%
\special{pn 20}%
\special{pa 2000 2000}%
\special{pa 2008 1964}%
\special{pa 2014 1930}%
\special{pa 2018 1896}%
\special{pa 2018 1862}%
\special{pa 2016 1832}%
\special{pa 2008 1804}%
\special{pa 1994 1780}%
\special{pa 1974 1760}%
\special{pa 1948 1742}%
\special{pa 1920 1724}%
\special{pa 1888 1710}%
\special{pa 1856 1696}%
\special{pa 1840 1690}%
\special{sp}%
%
\special{pn 8}%
\special{pa 1540 1540}%
\special{pa 1562 1516}%
\special{pa 1584 1492}%
\special{pa 1602 1466}%
\special{pa 1616 1438}%
\special{pa 1628 1410}%
\special{pa 1636 1378}%
\special{pa 1640 1346}%
\special{pa 1642 1314}%
\special{pa 1642 1280}%
\special{pa 1640 1246}%
\special{pa 1636 1212}%
\special{pa 1634 1176}%
\special{pa 1634 1144}%
\special{pa 1634 1110}%
\special{pa 1640 1080}%
\special{pa 1648 1052}%
\special{pa 1664 1026}%
\special{pa 1686 1002}%
\special{pa 1712 984}%
\special{pa 1744 968}%
\special{pa 1778 956}%
\special{pa 1812 948}%
\special{pa 1848 946}%
\special{pa 1882 948}%
\special{pa 1916 954}%
\special{pa 1944 968}%
\special{pa 1970 986}%
\special{pa 1994 1006}%
\special{pa 2016 1030}%
\special{pa 2038 1052}%
\special{pa 2062 1076}%
\special{pa 2086 1096}%
\special{pa 2116 1114}%
\special{pa 2148 1126}%
\special{pa 2184 1134}%
\special{pa 2220 1142}%
\special{pa 2254 1150}%
\special{pa 2288 1158}%
\special{pa 2316 1168}%
\special{pa 2338 1184}%
\special{pa 2352 1204}%
\special{pa 2358 1232}%
\special{pa 2360 1262}%
\special{pa 2358 1296}%
\special{pa 2356 1332}%
\special{pa 2358 1366}%
\special{pa 2366 1396}%
\special{pa 2380 1424}%
\special{pa 2398 1448}%
\special{pa 2422 1472}%
\special{pa 2450 1492}%
\special{pa 2478 1512}%
\special{pa 2490 1520}%
\special{sp}%
%
\special{pn 8}%
\special{pa 1260 1260}%
\special{pa 1272 1230}%
\special{pa 1284 1202}%
\special{pa 1296 1172}%
\special{pa 1308 1142}%
\special{pa 1320 1112}%
\special{pa 1332 1084}%
\special{pa 1346 1054}%
\special{pa 1360 1024}%
\special{pa 1374 996}%
\special{pa 1388 966}%
\special{pa 1402 938}%
\special{pa 1418 908}%
\special{pa 1436 880}%
\special{pa 1452 852}%
\special{pa 1470 826}%
\special{pa 1490 800}%
\special{pa 1510 774}%
\special{pa 1532 750}%
\special{pa 1554 728}%
\special{pa 1578 706}%
\special{pa 1602 686}%
\special{pa 1628 668}%
\special{pa 1656 652}%
\special{pa 1684 638}%
\special{pa 1714 624}%
\special{pa 1744 614}%
\special{pa 1774 604}%
\special{pa 1806 596}%
\special{pa 1838 590}%
\special{pa 1872 586}%
\special{pa 1904 582}%
\special{pa 1938 580}%
\special{pa 1972 580}%
\special{pa 2006 582}%
\special{pa 2038 584}%
\special{pa 2072 590}%
\special{pa 2106 596}%
\special{pa 2138 602}%
\special{pa 2170 610}%
\special{pa 2202 620}%
\special{pa 2232 632}%
\special{pa 2262 644}%
\special{pa 2292 658}%
\special{pa 2320 674}%
\special{pa 2348 690}%
\special{pa 2374 708}%
\special{pa 2400 728}%
\special{pa 2426 748}%
\special{pa 2450 770}%
\special{pa 2472 792}%
\special{pa 2494 816}%
\special{pa 2514 840}%
\special{pa 2534 864}%
\special{pa 2554 890}%
\special{pa 2574 914}%
\special{pa 2594 940}%
\special{pa 2614 964}%
\special{pa 2634 990}%
\special{pa 2654 1014}%
\special{pa 2674 1040}%
\special{pa 2694 1064}%
\special{pa 2714 1090}%
\special{pa 2734 1116}%
\special{pa 2754 1140}%
\special{pa 2772 1166}%
\special{pa 2792 1190}%
\special{pa 2800 1200}%
\special{sp}%
\put(21.0000,-10.8000){\makebox(0,0)[lb]{$f c_i'''$}}%
\put(28.1000,-7.0000){\makebox(0,0)[lb]{$U_1$}}%
%
\special{pn 8}%
\special{pa 3400 1200}%
\special{pa 4000 1200}%
\special{fp}%
%
\special{pn 8}%
\special{pa 3400 1310}%
\special{pa 4000 1310}%
\special{fp}%
%
\special{pn 8}%
\special{pa 4060 1250}%
\special{pa 3910 1370}%
\special{fp}%
%
\special{pn 8}%
\special{pa 4060 1250}%
\special{pa 3910 1130}%
\special{fp}%
%
\special{pn 8}%
\special{ar 5500 500 1500 1000  0.5300153 2.6183143}%
%
\special{pn 8}%
\special{sh 0}%
\special{ar 5880 1470 28 28  0.0000000 6.2831853}%
%
\special{pn 8}%
\special{sh 0.600}%
\special{ar 4520 1250 28 28  0.0000000 6.2831853}%
%
\special{pn 8}%
\special{sh 0.600}%
\special{ar 6500 1250 28 28  0.0000000 6.2831853}%
%
\special{pn 8}%
\special{sh 0.600}%
\special{ar 4850 1400 28 28  0.0000000 6.2831853}%
%
\special{pn 8}%
\special{sh 0.600}%
\special{ar 6220 1380 28 28  0.0000000 6.2831853}%
%
\special{pn 8}%
\special{sh 0}%
\special{ar 5210 1480 28 28  0.0000000 6.2831853}%
\put(52.2000,-14.2000){\makebox(0,0)[lb]{$\eta_1^\infty$}}%
\put(48.4000,-14.3000){\makebox(0,0)[rt]{$\widehat{f_1}^n\eta_1$}}%
\put(62.3000,-14.1000){\makebox(0,0)[lt]{$\widehat{f_1}^n\zeta_1$}}%
\put(55.0000,-10.0000){\makebox(0,0){$\widehat{U_1}$}}%
\put(58.4000,-14.3000){\makebox(0,0)[lb]{$\zeta_1^\infty$}}%
\put(16.3000,-5.5000){\makebox(0,0)[lb]{$c_i'''$}}%
\end{picture}%
\caption{No topological chain $\{c_i'''\}$}
\label{fig14}
\end{center}
\end{figure}
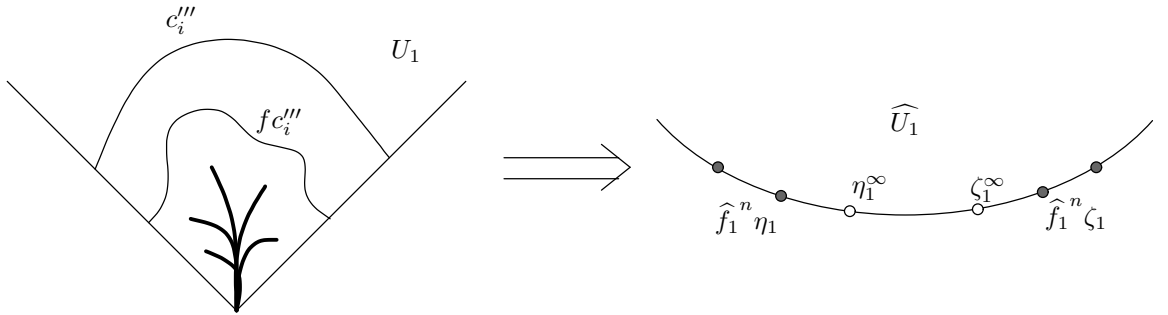

Let us show the last part of the theorem. Assume $\xi$
is an attractor of $\hat f\vert_{\PP(U)}$. Choose an extendable
topological chain $\{c_i\}$ representing $\xi$. Then as before we can assume
$fU(c_i)\subset U(c_i)$ and $U(c_i)\cap{\rm Fix}(f)=\emptyset$
for any big $i$. Fix one such $i$ and let $c=c_i$. 
Let
$U_1=U\setminus\cap_{n\geq 1}f^nU(c)$. See Figure \ref{fig14}.
Our purpose is to show that $U_1=U$. Notice that this implies that
$\xi$ is an attractor of $\hat f$.
Denote the two end points of $c$ in $\PP(U_1)$ by $\eta_1$ and 
$\zeta_1$ and let $\eta_1^\infty=\lim \hat f_1^n\eta_1$
and
$\zeta_1^\infty=\lim \hat f_1^n\zeta_1$,
where $\hat f_1$ is the homeomorphism of $\hat U_1$ induced by
$f$. 
We have $\zeta_1^\infty=\eta_1^\infty$, for otherwise
the same argument as before yields a contradiction.
 Take an extendable
topological chain $\{c'''_i\}$ representing this prime end in
$\PP(U_1)$.
It is also a topological chain for $U$ and we have
$$
U\setminus U(c'''_i)=U_1\setminus U_1(c'''_i).
$$
Since $\cap_iU(c'''_i)=\cap_iU_1(c'''_i)=\emptyset$ by Lemma \ref{0a}, 
this shows $U_1=U$, as is required.
\qed

\section{Minimal continuum}

Let $f$ be an orientation preserving homeomorphism of the 2-sphere $S^2$ 
which
has a continuum $X$ as a minimal set. Recall that
a connected
component $U$ of $S^2\setminus X$ is called an invariant domain
if $fU=U$.
The purpose of this section is to prove Theorem \ref{0}.
We begin with the following lemma.

\begin{lemma} \label{9}
The Carath\'eodory rotation number of an invariant domain $U$ is nonzero.
\end{lemma}

Before the proof, let us mention that Example \ref{example1} shows
the necessity for the minimality assumption and that Example
\ref{example2}
shows that Lemma \ref{9} does not hold for surfaces of nonzero genus.

\bigskip
{\bf Proof of Lemma \ref{9}}.
Denote by $\hat f$ the homeomorphism that $f$ induces on $\hat U$.
Assume for contradiction that the rotation number of 
$\hat f\vert_{\PP(U)}$ is 0. Then the conclusion of Theorem \ref{5}
holds.
Let $\alpha$ and $\omega$ be adjacent repelling and attracting
fixed points on $\PP(U)$ and choose an interval $(\alpha,\omega)$
in $\PP(U)$ so that $(\alpha,\omega)\cap{\rm Fix}(\hat f)=\emptyset$.
By Lemma \ref{4} there is
a prime end $\xi\in(\alpha,\omega)$ belonging
to the set $\EE(U)$ of the extendable prime ends near $\omega$.
Then one can choose an extendable curve $\hat\gamma$ joining $\xi$
and $\hat f\xi$ such that $\gamma=\hat\gamma\cap U$ is contained
in an open fundamental domain $F$ of $\hat f$. (Recall that $\omega$ is
an attractor of the homeomorphism $\hat f$.) See Figure \ref{fig15}. 
\begin{figure}[htbp]
\begin{center}
\unitlength 0.1in
\begin{picture}( 41.0800, 26.0800)(  4.2000,-30.8000)
\put(35.1000,-23.2000){\makebox(0,0)[lb]{$\gamma$}}%
%
\special{pn 8}%
\special{pa 500 500}%
\special{pa 4500 500}%
\special{pa 4500 3000}%
\special{pa 500 3000}%
\special{pa 500 500}%
\special{fp}%
%
\special{pn 8}%
\special{pa 4500 500}%
\special{pa 2500 500}%
\special{fp}%
\special{sh 1}%
\special{pa 2500 500}%
\special{pa 2568 520}%
\special{pa 2554 500}%
\special{pa 2568 480}%
\special{pa 2500 500}%
\special{fp}%
%
\special{pn 8}%
\special{pa 500 3000}%
\special{pa 2500 3000}%
\special{fp}%
\special{sh 1}%
\special{pa 2500 3000}%
\special{pa 2434 2980}%
\special{pa 2448 3000}%
\special{pa 2434 3020}%
\special{pa 2500 3000}%
\special{fp}%
%
\special{pn 8}%
\special{pa 500 3000}%
\special{pa 500 1750}%
\special{fp}%
\special{sh 1}%
\special{pa 500 1750}%
\special{pa 480 1818}%

\special{pa 500 1804}%
\special{pa 520 1818}%
\special{pa 500 1750}%
\special{fp}%
%
\special{pn 8}%
\special{ar 1620 3000 300 1250  3.1415927 6.2831853}%
%
\special{pn 8}%
\special{ar 3550 3250 270 930  3.4256939 5.9990840}%
%
\special{pn 8}%
\special{ar 2930 3250 350 1120  3.3608666 6.0639114}%
%
\special{pn 8}%
\special{ar 2250 3250 340 1300  3.3304898 6.0942881}%
%
\special{pn 8}%
\special{ar 3940 3250 150 690  3.5352909 5.8894871}%
%
\special{pn 8}%
\special{pa 3280 3000}%
\special{pa 3280 2968}%
\special{pa 3280 2936}%
\special{pa 3280 2904}%
\special{pa 3280 2872}%
\special{pa 3280 2840}%
\special{pa 3280 2808}%
\special{pa 3280 2776}%
\special{pa 3280 2744}%
\special{pa 3280 2712}%
\special{pa 3280 2680}%
\special{pa 3280 2648}%
\special{pa 3280 2616}%
\special{pa 3280 2584}%
\special{pa 3278 2552}%
\special{pa 3278 2520}%
\special{pa 3278 2488}%
\special{pa 3280 2456}%
\special{pa 3280 2424}%
\special{pa 3280 2392}%
\special{pa 3282 2360}%
\special{pa 3282 2328}%
\special{pa 3284 2296}%
\special{pa 3286 2264}%
\special{pa 3290 2232}%
\special{pa 3294 2200}%
\special{pa 3300 2168}%
\special{pa 3308 2138}%
\special{pa 3316 2106}%
\special{pa 3326 2076}%
\special{pa 3338 2046}%
\special{pa 3350 2016}%
\special{pa 3366 1988}%
\special{pa 3382 1960}%
\special{pa 3400 1934}%
\special{pa 3420 1908}%
\special{pa 3440 1882}%
\special{pa 3462 1858}%
\special{pa 3486 1836}%
\special{pa 3510 1814}%
\special{pa 3534 1794}%
\special{pa 3560 1776}%
\special{pa 3588 1758}%
\special{pa 3614 1742}%
\special{pa 3644 1726}%
\special{pa 3672 1712}%
\special{pa 3702 1698}%
\special{pa 3732 1686}%
\special{pa 3762 1674}%
\special{pa 3792 1664}%
\special{pa 3822 1654}%
\special{pa 3854 1646}%
\special{pa 3884 1638}%
\special{pa 3916 1630}%
\special{pa 3948 1624}%
\special{pa 3980 1618}%
\special{pa 4012 1614}%
\special{pa 4044 1610}%
\special{pa 4076 1608}%
\special{pa 4108 1604}%
\special{pa 4140 1602}%
\special{pa 4172 1602}%
\special{pa 4204 1600}%
\special{pa 4236 1600}%
\special{pa 4268 1600}%
\special{pa 4300 1600}%
\special{pa 4332 1602}%
\special{pa 4364 1602}%
\special{pa 4396 1604}%
\special{pa 4428 1606}%
\special{pa 4460 1608}%
\special{pa 4492 1610}%
\special{pa 4500 1610}%
\special{sp 0.070}%
\put(40.3000,-19.5000){\makebox(0,0)[lb]{$F$}}%
\put(28.3000,-21.2000){\makebox(0,0)[lb]{$f^{-1}\gamma$}}%
\put(44.2000,-30.8000){\makebox(0,0)[lt]{$\omega$}}%
\put(37.6000,-30.5000){\makebox(0,0)[lt]{$\widehat{f}\xi$}}%
\put(32.5000,-30.4000){\makebox(0,0)[lt]{$\xi$}}%
%
\special{pn 8}%
\special{sh 0.600}%
\special{ar 4500 3000 28 28  0.0000000 6.2831853}%
%
\special{pn 8}%
\special{sh 0.600}%
\special{ar 500 500 28 28  0.0000000 6.2831853}%
%
\special{pn 8}%
\special{sh 0}%
\special{ar 500 3000 28 28  0.0000000 6.2831853}%
%
\special{pn 8}%
\special{pa 4500 500}%
\special{pa 4500 1770}%
\special{fp}%
\special{sh 1}%
\special{pa 4500 1770}%
\special{pa 4520 1704}%
\special{pa 4500 1718}%
\special{pa 4480 1704}%
\special{pa 4500 1770}%
\special{fp}%
%
\special{pn 8}%
\special{pa 3810 3000}%
\special{pa 3810 2968}%
\special{pa 3810 2936}%
\special{pa 3810 2904}%
\special{pa 3810 2872}%
\special{pa 3810 2840}%
\special{pa 3810 2808}%
\special{pa 3810 2776}%
\special{pa 3810 2744}%
\special{pa 3810 2712}%
\special{pa 3810 2680}%
\special{pa 3810 2648}%
\special{pa 3810 2616}%
\special{pa 3810 2584}%
\special{pa 3812 2552}%
\special{pa 3812 2520}%
\special{pa 3814 2486}%
\special{pa 3818 2454}%
\special{pa 3822 2422}%
\special{pa 3830 2390}%
\special{pa 3838 2360}%
\special{pa 3850 2330}%
\special{pa 3866 2302}%
\special{pa 3882 2276}%
\special{pa 3902 2250}%
\special{pa 3926 2228}%
\special{pa 3950 2206}%
\special{pa 3976 2186}%
\special{pa 4004 2170}%
\special{pa 4032 2156}%
\special{pa 4062 2144}%
\special{pa 4094 2134}%
\special{pa 4124 2126}%
\special{pa 4156 2120}%
\special{pa 4188 2116}%
\special{pa 4220 2112}%
\special{pa 4252 2112}%
\special{pa 4284 2110}%
\special{pa 4316 2110}%
\special{pa 4348 2110}%
\special{pa 4380 2110}%
\special{pa 4412 2110}%
\special{pa 4444 2110}%
\special{pa 4476 2110}%
\special{pa 4500 2110}%
\special{sp 0.070}%
\put(4.2000,-30.8000){\makebox(0,0)[lt]{$\alpha$}}%
%
\special{pn 8}%
\special{ar 3280 3000 650 650  4.4937200 4.9759529}%
%
\special{pn 8}%
\special{pa 2920 2610}%
\special{pa 3240 2440}%
\special{fp}%
\special{sh 1}%
\special{pa 3240 2440}%
\special{pa 3172 2454}%
\special{pa 3194 2466}%
\special{pa 3192 2490}%
\special{pa 3240 2440}%
\special{fp}%
\put(29.2000,-26.2000){\makebox(0,0)[rt]{$\hat V$}}%
%
\special{pn 8}%
\special{sh 0}%
\special{ar 4500 500 28 28  0.0000000 6.2831853}%
\end{picture}%
\caption{$\widehat{U}$}
\label{fig15}
\end{center}
\end{figure}
Notice that the natural map $\Phi:\EE(U)\to X$ is equivariant,
$f\circ\Phi=\Phi\circ\hat f$.
Therefore 
the closure $\overline\gamma$ of the curve $\gamma$
in $S^2$ joins a point, say $p$, with $fp$.
Notice that $p\in X$.
The cross cuts $f^n\gamma$ in $U$
($n\in{\mathbb Z}$) are
mutually disjoint and its closure 
$f^n(\overline\gamma)$ joins a point $f^n(p)$ with $f^{n+1}(p)$.

Since $X$ is minimal and $p\in X$, there is $n>0$ such that
$f^np$ is arbitrarily near $p$. 
Consider a small disc $B$ centered at $p$ such
that $B\cap fB=\emptyset$. The connected
component
of $f^{-1}\overline\gamma\cup\overline\gamma$ that contains the point
$p$ divides $B$ into two domains. One of it $V$, corresponding to $\hat V$ in Figure \ref{fig15}, is contained in $U$ (if we choose $B$ small enough) and
the point $f^np$ can be chosen from the component
of $B\setminus(f^{-1}\overline \gamma\cup \overline\gamma)$
adjacent to $V$.
Choose a small arc $\delta'$ in $B$ joining $p$
with $f^np$ which does not intersect 
$f^{-1}\overline\gamma\cup\overline\gamma$ except at $p$.
Notice
that $f\delta'\cap\delta'=\emptyset$. See Figure \ref{fig16}.
\begin{figure}[htbp]
\begin{center}
\unitlength 0.1in
\begin{picture}( 60.0400, 30.0000)(  1.4000,-35.0000)
%
\special{pn 8}%
\special{ar 2000 2000 1500 1500  0.0000000 6.2831853}%
%
\special{pn 8}%
\special{sh 0.600}%
\special{ar 2000 2250 28 28  0.0000000 6.2831853}%
%
\special{pn 8}%
\special{sh 0.600}%
\special{ar 2000 2000 28 28  0.0000000 6.2831853}%
%
\special{pn 8}%
\special{sh 0.600}%
\special{ar 5278 2238 28 28  0.0000000 6.2831853}%
%
\special{pn 8}%
\special{sh 0.600}%
\special{ar 5278 1988 28 28  0.0000000 6.2831853}%
%
\special{pn 8}%
\special{pa 2000 2000}%
\special{pa 1980 1976}%
\special{pa 1960 1950}%
\special{pa 1938 1926}%
\special{pa 1916 1902}%
\special{pa 1894 1880}%
\special{pa 1872 1858}%
\special{pa 1848 1836}%
\special{pa 1824 1816}%
\special{pa 1798 1796}%
\special{pa 1772 1778}%
\special{pa 1744 1762}%
\special{pa 1716 1746}%
\special{pa 1686 1732}%
\special{pa 1656 1720}%
\special{pa 1626 1708}%
\special{pa 1594 1698}%
\special{pa 1562 1690}%
\special{pa 1530 1684}%
\special{pa 1498 1678}%
\special{pa 1466 1674}%
\special{pa 1434 1670}%
\special{pa 1400 1670}%
\special{pa 1368 1670}%
\special{pa 1336 1672}%
\special{pa 1304 1674}%
\special{pa 1272 1678}%
\special{pa 1240 1684}%
\special{pa 1208 1690}%
\special{pa 1176 1698}%
\special{pa 1144 1706}%
\special{pa 1112 1714}%
\special{pa 1082 1724}%
\special{pa 1052 1736}%
\special{pa 1022 1746}%
\special{pa 992 1758}%
\special{pa 962 1772}%
\special{pa 934 1784}%
\special{pa 906 1798}%
\special{pa 878 1812}%
\special{pa 850 1826}%
\special{pa 824 1840}%
\special{pa 798 1856}%
\special{pa 772 1870}%
\special{pa 746 1886}%
\special{pa 722 1904}%
\special{pa 698 1920}%
\special{pa 674 1938}%
\special{pa 650 1958}%
\special{pa 626 1976}%
\special{pa 604 1996}%
\special{pa 582 2018}%
\special{pa 560 2038}%
\special{pa 538 2062}%
\special{pa 518 2084}%
\special{pa 496 2108}%
\special{pa 476 2134}%
\special{pa 456 2160}%
\special{pa 436 2188}%
\special{pa 418 2216}%
\special{pa 398 2246}%
\special{pa 380 2276}%
\special{pa 362 2308}%
\special{pa 344 2340}%
\special{pa 326 2374}%
\special{pa 308 2410}%
\special{pa 292 2446}%
\special{pa 274 2484}%
\special{pa 258 2522}%
\special{pa 246 2558}%
\special{pa 238 2590}%
\special{pa 238 2612}%
\special{pa 246 2626}%
\special{pa 264 2628}%
\special{pa 288 2620}%
\special{pa 318 2606}%
\special{pa 350 2586}%
\special{pa 386 2560}%
\special{pa 418 2532}%
\special{pa 450 2502}%
\special{pa 476 2472}%
\special{pa 498 2444}%
\special{pa 516 2414}%
\special{pa 532 2386}%
\special{pa 546 2358}%
\special{pa 560 2332}%
\special{pa 572 2306}%
\special{pa 586 2280}%
\special{pa 600 2254}%
\special{pa 618 2230}%
\special{pa 636 2208}%
\special{pa 658 2186}%
\special{pa 682 2164}%
\special{pa 706 2142}%
\special{pa 732 2122}%
\special{pa 758 2102}%
\special{pa 786 2084}%
\special{pa 814 2066}%
\special{pa 842 2048}%
\special{pa 870 2030}%
\special{pa 896 2014}%
\special{pa 922 1996}%
\special{pa 948 1980}%
\special{pa 974 1966}%
\special{pa 1000 1952}%
\special{pa 1028 1938}%
\special{pa 1058 1926}%
\special{pa 1090 1914}%
\special{pa 1126 1904}%
\special{pa 1164 1896}%
\special{pa 1204 1890}%
\special{pa 1242 1886}%
\special{pa 1276 1886}%
\special{pa 1300 1892}%
\special{pa 1312 1904}%
\special{pa 1310 1922}%
\special{pa 1296 1944}%
\special{pa 1274 1970}%
\special{pa 1244 1998}%
\special{pa 1210 2024}%
\special{pa 1176 2048}%
\special{pa 1140 2068}%
\special{pa 1108 2086}%
\special{pa 1076 2100}%
\special{pa 1046 2112}%
\special{pa 1016 2122}%
\special{pa 988 2132}%
\special{pa 960 2142}%
\special{pa 934 2154}%
\special{pa 910 2166}%
\special{pa 886 2180}%
\special{pa 862 2196}%
\special{pa 838 2214}%
\special{pa 816 2234}%
\special{pa 796 2256}%
\special{pa 774 2278}%
\special{pa 754 2302}%
\special{pa 734 2328}%
\special{pa 714 2354}%
\special{pa 694 2382}%
\special{pa 674 2408}%
\special{pa 656 2438}%
\special{pa 636 2466}%
\special{pa 618 2494}%
\special{pa 600 2524}%
\special{pa 580 2552}%
\special{pa 560 2580}%
\special{pa 542 2608}%
\special{pa 522 2636}%
\special{pa 502 2662}%
\special{pa 480 2688}%
\special{pa 460 2712}%
\special{pa 438 2736}%
\special{pa 416 2760}%
\special{pa 394 2782}%
\special{pa 370 2802}%
\special{pa 344 2822}%
\special{pa 320 2840}%
\special{pa 294 2858}%
\special{pa 266 2874}%
\special{pa 238 2888}%
\special{pa 208 2900}%
\special{pa 178 2914}%
\special{pa 160 2920}%
\special{sp}%
%
\special{pn 8}%
\special{pa 2000 2000}%
\special{pa 2000 2250}%
\special{fp}%
\put(20.0000,-11.5000){\makebox(0,0){$V$}}%
\put(20.0000,-18.2000){\makebox(0,0){$p$}}%
\put(19.8000,-22.1000){\makebox(0,0)[rb]{$\delta'$}}%
\put(20.4000,-23.0000){\makebox(0,0)[lb]{$f^np$}}%
\put(1.4000,-21.6000){\makebox(0,0)[lb]{$f^{-1}\overline{\gamma}$}}%
%
\special{pn 8}%
\special{ar 5268 2118 878 878  0.0000000 6.2831853}%
\put(53.2700,-19.5700){\makebox(0,0)[lb]{$fp$}}%
\put(53.4000,-21.6000){\makebox(0,0)[lb]{$f\delta'$}}%
\put(51.9000,-12.0000){\makebox(0,0)[lb]{$fB$}}%
\put(39.9000,-13.9000){\makebox(0,0)[lb]{$\overline{\gamma}$}}%
\put(28.1000,-6.9000){\makebox(0,0)[lb]{$B$}}%
%
\special{pn 8}%
\special{pa 2000 2000}%
\special{pa 2022 1978}%
\special{pa 2046 1954}%
\special{pa 2068 1932}%
\special{pa 2090 1908}%
\special{pa 2112 1886}%
\special{pa 2136 1864}%
\special{pa 2160 1842}%
\special{pa 2182 1820}%
\special{pa 2206 1800}%
\special{pa 2232 1778}%
\special{pa 2256 1758}%
\special{pa 2280 1738}%
\special{pa 2306 1718}%
\special{pa 2332 1700}%
\special{pa 2358 1680}%
\special{pa 2384 1662}%
\special{pa 2410 1644}%
\special{pa 2438 1626}%
\special{pa 2464 1608}%
\special{pa 2492 1592}%
\special{pa 2520 1576}%
\special{pa 2548 1560}%
\special{pa 2576 1544}%
\special{pa 2604 1530}%
\special{pa 2632 1516}%
\special{pa 2662 1502}%
\special{pa 2692 1490}%
\special{pa 2722 1478}%
\special{pa 2752 1468}%
\special{pa 2782 1456}%
\special{pa 2812 1448}%
\special{pa 2844 1440}%
\special{pa 2874 1432}%
\special{pa 2906 1424}%
\special{pa 2938 1418}%
\special{pa 2970 1412}%
\special{pa 3002 1406}%
\special{pa 3034 1402}%
\special{pa 3066 1398}%
\special{pa 3098 1394}%
\special{pa 3130 1390}%
\special{pa 3162 1386}%
\special{pa 3194 1382}%
\special{pa 3226 1380}%
\special{pa 3258 1376}%
\special{pa 3288 1372}%
\special{pa 3320 1370}%
\special{pa 3352 1366}%
\special{pa 3382 1362}%
\special{pa 3414 1358}%
\special{pa 3444 1354}%
\special{pa 3474 1350}%
\special{pa 3506 1348}%
\special{pa 3536 1344}%
\special{pa 3566 1340}%
\special{pa 3598 1338}%
\special{pa 3628 1336}%
\special{pa 3660 1336}%
\special{pa 3692 1338}%
\special{pa 3724 1340}%
\special{pa 3756 1342}%
\special{pa 3790 1348}%
\special{pa 3824 1356}%
\special{pa 3860 1364}%
\special{pa 3896 1376}%
\special{pa 3932 1390}%
\special{pa 3968 1404}%
\special{pa 4004 1422}%
\special{pa 4036 1440}%
\special{pa 4066 1462}%
\special{pa 4090 1482}%
\special{pa 4110 1504}%
\special{pa 4124 1528}%
\special{pa 4130 1552}%
\special{pa 4128 1576}%
\special{pa 4118 1598}%
\special{pa 4102 1622}%
\special{pa 4080 1644}%
\special{pa 4054 1666}%
\special{pa 4022 1688}%
\special{pa 3990 1706}%
\special{pa 3954 1722}%
\special{pa 3920 1738}%
\special{pa 3886 1748}%
\special{pa 3852 1758}%
\special{pa 3818 1766}%
\special{pa 3786 1770}%
\special{pa 3754 1774}%
\special{pa 3722 1776}%
\special{pa 3692 1776}%
\special{pa 3660 1776}%
\special{pa 3630 1776}%
\special{pa 3598 1774}%
\special{pa 3568 1774}%
\special{pa 3536 1772}%
\special{pa 3504 1770}%
\special{pa 3472 1770}%
\special{pa 3440 1772}%
\special{pa 3406 1772}%
\special{pa 3374 1774}%
\special{pa 3340 1778}%
\special{pa 3306 1782}%
\special{pa 3274 1788}%
\special{pa 3240 1796}%
\special{pa 3208 1804}%
\special{pa 3178 1812}%
\special{pa 3146 1824}%
\special{pa 3116 1836}%
\special{pa 3088 1850}%
\special{pa 3060 1866}%
\special{pa 3034 1882}%
\special{pa 3010 1900}%
\special{pa 2988 1922}%
\special{pa 2966 1944}%
\special{pa 2946 1966}%
\special{pa 2926 1992}%
\special{pa 2908 2018}%
\special{pa 2892 2044}%
\special{pa 2876 2072}%
\special{pa 2860 2100}%
\special{pa 2844 2128}%
\special{pa 2828 2158}%
\special{pa 2812 2186}%
\special{pa 2798 2216}%
\special{pa 2782 2244}%
\special{pa 2764 2274}%
\special{pa 2748 2302}%
\special{pa 2728 2328}%
\special{pa 2710 2356}%
\special{pa 2690 2380}%
\special{pa 2668 2406}%
\special{pa 2646 2428}%
\special{pa 2622 2452}%
\special{pa 2598 2474}%
\special{pa 2574 2494}%
\special{pa 2548 2512}%
\special{pa 2520 2532}%
\special{pa 2494 2548}%
\special{pa 2466 2564}%
\special{pa 2436 2578}%
\special{pa 2408 2592}%
\special{pa 2378 2604}%
\special{pa 2348 2614}%
\special{pa 2316 2624}%
\special{pa 2286 2632}%
\special{pa 2254 2638}%
\special{pa 2222 2644}%
\special{pa 2190 2648}%
\special{pa 2158 2652}%
\special{pa 2126 2654}%
\special{pa 2092 2654}%
\special{pa 2060 2654}%
\special{pa 2028 2652}%
\special{pa 1996 2648}%
\special{pa 1962 2644}%
\special{pa 1930 2638}%
\special{pa 1898 2630}%
\special{pa 1866 2622}%
\special{pa 1836 2612}%
\special{pa 1804 2602}%
\special{pa 1774 2590}%
\special{pa 1744 2576}%
\special{pa 1716 2562}%
\special{pa 1688 2546}%
\special{pa 1660 2528}%
\special{pa 1634 2510}%
\special{pa 1610 2490}%
\special{pa 1586 2468}%
\special{pa 1562 2446}%
\special{pa 1540 2424}%
\special{pa 1518 2402}%
\special{pa 1494 2380}%
\special{pa 1470 2358}%
\special{pa 1446 2338}%
\special{pa 1418 2320}%
\special{pa 1390 2304}%
\special{pa 1360 2290}%
\special{pa 1330 2278}%
\special{pa 1298 2268}%
\special{pa 1266 2262}%
\special{pa 1234 2258}%
\special{pa 1202 2258}%
\special{pa 1170 2260}%
\special{pa 1140 2266}%
\special{pa 1110 2274}%
\special{pa 1080 2286}%
\special{pa 1052 2300}%
\special{pa 1024 2316}%
\special{pa 998 2332}%
\special{pa 970 2352}%
\special{pa 944 2372}%
\special{pa 918 2394}%
\special{pa 894 2416}%
\special{pa 868 2438}%
\special{pa 844 2460}%
\special{pa 818 2484}%
\special{pa 794 2506}%
\special{pa 770 2528}%
\special{pa 746 2550}%
\special{pa 722 2572}%
\special{pa 698 2594}%
\special{pa 676 2616}%
\special{pa 656 2640}%
\special{pa 638 2666}%
\special{pa 620 2692}%
\special{pa 604 2720}%
\special{pa 590 2750}%
\special{pa 578 2782}%
\special{pa 570 2816}%
\special{pa 562 2854}%
\special{pa 560 2890}%
\special{pa 560 2928}%
\special{pa 564 2962}%
\special{pa 574 2994}%
\special{pa 588 3020}%
\special{pa 608 3038}%
\special{pa 634 3050}%
\special{pa 664 3054}%
\special{pa 698 3052}%
\special{pa 732 3044}%
\special{pa 766 3032}%
\special{pa 798 3016}%
\special{pa 824 2998}%
\special{pa 850 2978}%
\special{pa 872 2956}%
\special{pa 894 2932}%
\special{pa 918 2910}%
\special{pa 940 2890}%
\special{pa 966 2870}%
\special{pa 992 2852}%
\special{pa 1020 2838}%
\special{pa 1052 2826}%
\special{pa 1082 2814}%
\special{pa 1114 2808}%
\special{pa 1146 2802}%
\special{pa 1178 2800}%
\special{pa 1210 2800}%
\special{pa 1242 2802}%
\special{pa 1274 2806}%
\special{pa 1304 2810}%
\special{pa 1336 2818}%
\special{pa 1368 2826}%
\special{pa 1398 2834}%
\special{pa 1430 2844}%
\special{pa 1462 2854}%
\special{pa 1492 2864}%
\special{pa 1524 2874}%
\special{pa 1554 2884}%
\special{pa 1586 2894}%
\special{pa 1618 2902}%
\special{pa 1648 2912}%
\special{pa 1680 2920}%
\special{pa 1712 2928}%
\special{pa 1742 2936}%
\special{pa 1774 2944}%
\special{pa 1806 2950}%
\special{pa 1836 2958}%
\special{pa 1868 2962}%
\special{pa 1900 2968}%
\special{pa 1932 2972}%
\special{pa 1962 2976}%
\special{pa 1994 2978}%
\special{pa 2026 2980}%
\special{pa 2058 2980}%
\special{pa 2088 2980}%
\special{pa 2120 2980}%
\special{pa 2152 2978}%
\special{pa 2184 2974}%
\special{pa 2216 2970}%
\special{pa 2246 2966}%
\special{pa 2278 2960}%
\special{pa 2310 2952}%
\special{pa 2342 2946}%
\special{pa 2372 2938}%
\special{pa 2404 2928}%
\special{pa 2434 2918}%
\special{pa 2464 2908}%
\special{pa 2496 2896}%
\special{pa 2526 2884}%
\special{pa 2556 2870}%
\special{pa 2586 2858}%
\special{pa 2616 2844}%
\special{pa 2644 2828}%
\special{pa 2674 2814}%
\special{pa 2702 2798}%
\special{pa 2730 2780}%
\special{pa 2760 2764}%
\special{pa 2786 2746}%
\special{pa 2814 2728}%
\special{pa 2842 2710}%
\special{pa 2868 2690}%
\special{pa 2894 2672}%
\special{pa 2920 2652}%
\special{pa 2946 2632}%
\special{pa 2972 2612}%
\special{pa 2996 2592}%
\special{pa 3022 2572}%
\special{pa 3046 2552}%
\special{pa 3072 2530}%
\special{pa 3096 2510}%
\special{pa 3122 2490}%
\special{pa 3146 2470}%
\special{pa 3172 2450}%
\special{pa 3196 2430}%
\special{pa 3222 2412}%
\special{pa 3248 2392}%
\special{pa 3274 2374}%
\special{pa 3298 2356}%
\special{pa 3326 2338}%
\special{pa 3352 2322}%
\special{pa 3378 2304}%
\special{pa 3406 2290}%
\special{pa 3432 2274}%
\special{pa 3460 2260}%
\special{pa 3488 2246}%
\special{pa 3518 2232}%
\special{pa 3546 2220}%
\special{pa 3576 2208}%
\special{pa 3606 2196}%
\special{pa 3636 2186}%
\special{pa 3666 2174}%
\special{pa 3696 2164}%
\special{pa 3726 2154}%
\special{pa 3758 2144}%
\special{pa 3788 2136}%
\special{pa 3820 2126}%
\special{pa 3850 2116}%
\special{pa 3882 2108}%
\special{pa 3912 2098}%
\special{pa 3944 2090}%
\special{pa 3976 2080}%
\special{pa 4006 2072}%
\special{pa 4038 2062}%
\special{pa 4068 2054}%
\special{pa 4100 2044}%
\special{pa 4132 2034}%
\special{pa 4162 2024}%
\special{pa 4194 2014}%
\special{pa 4224 2006}%
\special{pa 4256 1996}%
\special{pa 4286 1986}%
\special{pa 4318 1976}%
\special{pa 4348 1966}%
\special{pa 4378 1956}%
\special{pa 4410 1948}%
\special{pa 4440 1938}%
\special{pa 4470 1928}%
\special{pa 4500 1918}%
\special{pa 4532 1910}%
\special{pa 4562 1900}%
\special{pa 4592 1892}%
\special{pa 4622 1884}%
\special{pa 4652 1874}%
\special{pa 4682 1866}%
\special{pa 4712 1858}%
\special{pa 4742 1850}%
\special{pa 4772 1842}%
\special{pa 4802 1836}%
\special{pa 4832 1828}%
\special{pa 4862 1822}%
\special{pa 4892 1816}%
\special{pa 4924 1812}%
\special{pa 4954 1808}%
\special{pa 4986 1804}%
\special{pa 5020 1802}%
\special{pa 5054 1800}%
\special{pa 5088 1800}%
\special{pa 5124 1802}%
\special{pa 5160 1804}%
\special{pa 5194 1810}%
\special{pa 5226 1822}%
\special{pa 5250 1838}%
\special{pa 5266 1862}%
\special{pa 5276 1892}%
\special{pa 5280 1926}%
\special{pa 5282 1960}%
\special{pa 5282 1992}%
\special{pa 5282 2024}%
\special{pa 5282 2054}%
\special{pa 5280 2086}%
\special{pa 5280 2118}%
\special{pa 5280 2150}%
\special{pa 5280 2182}%
\special{pa 5280 2214}%
\special{pa 5280 2240}%
\special{sp}%
\end{picture}%
\caption{}
\label{fig16}
\end{center}
\end{figure}

\begin{figure}[htbp]
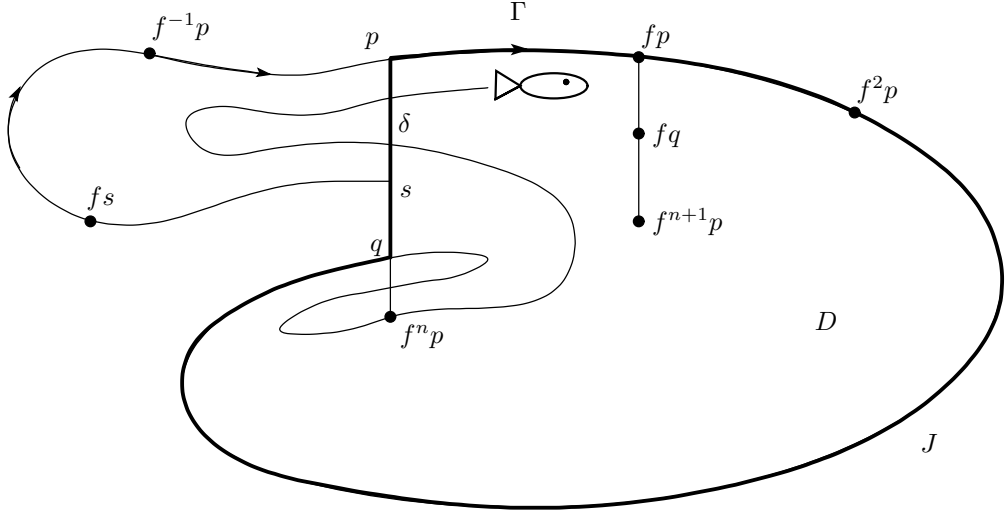

\begin{center}
\unitlength 0.1in
%
\caption{the curve $J$}
\label{fig17}
\end{center}
\end{figure}

Consider a long simple curve $\Gamma_+=\cup_{n\geq0}f^n\overline\gamma$.
Let $q$ be the first point of intersection 
of $\Gamma_+\setminus\{p\}$ with $\delta'$ (possibly $q=f^np$)
and let $\delta$ be the subarc of
$\delta'$ joining $p$ and $q$. Notice that $q$ is not from 
$\overline\gamma$ since $\delta'\cap\overline\gamma=\{p\}$.
The tiny arc $\delta$ together with the
subarc $\Gamma_+^0$ of $\Gamma_+$ that joins $p$ and $q$ forms a Jordan
curve
$J$.
See Figure \ref{fig17}.

Let $D$ be the connected component of $S^2\setminus J$ which
contains $fq$. Then the half open arc $f\delta'\setminus\{fp\}$
cannot intersect $J$ since $q$ is the first intersection point.
Thus $f\delta'\setminus\{fp\}$ and {\em in particular}
its end point $f^{n+1}p$ is contained in
$D$. 

We also have $f^{-1}\gamma\cap \overline D=\emptyset$.
In fact
$f^{-1}$ is an orientation preserving
 homeomorphism mapping a neighbourhood
of $fp$ to a neighbourhood of $p$. So the cyclic order of the
three curves $\overline\gamma$, $f\delta$, $f\overline\gamma$ emanating from the
point $fp$ is the same as the cyclic order of the curves
$f^{-1}\overline\gamma$, $\delta$, $\overline\gamma$ emanating from $p$. That is,
the curve $f^{-1}\overline\gamma$ tends towards outside of $D$, and thus
$f^{-1}\gamma\cap \overline D=\emptyset$.

Another long curve $\Gamma_-=\cup_{n<0}f^n\overline\gamma$
must go arbitrarily near to the point $f^{n+1}p$ which is in $D$, and therefore
must
intersect $\delta$.
Let
$s$ be the first intersection point of $\Gamma_-\setminus\{p\}$ with $\delta$.
Then an open arc $\Gamma_-^0$ in $\Gamma_-$ with end points $p$ and $s$ cannot
intersect $J$ and therefore $\Gamma_-^0\cap D=\emptyset$. 
By the construction of $\delta'$, $s$ is not from
$f^{-1}\overline\gamma$
and thus  $fs\in\Gamma_-^0$.
On the other hand $fs$ lies on $f\delta$ and therefore
belongs to $D$. A contradiction.
\qed
 
\bigskip
A closed disc $D$ in $S^2$ is called {\em adapted} if 
$\partial D\cap{\rm Fix}(f)=\emptyset$ 
and $D\cup fD\neq S^2$. Given an adapted disc $D$,
choosing the point of infinity in $S^2\setminus(D\cup fD)$,
one may consider $D\cup fD$ to be contained in $\R^2$.
Then the degree of the map
$$
id-f:\partial D\to\R^2\setminus\{0\}
$$
is called the {\em index} of $f$ w.\ r.\ t.\ $D$
and is denoted by $\ii D$.
An application of the Lefschetz index theorem yields
the following lemma.

\begin{lemma} \label{10}
Let $D_1,\cdots,D_r$ be mutually disjoint adapted discs
such that there is no fixed point of $f$ in the complement
of $\cup_{j=1}^rD_j$. Then we have
$$
\sum_{j=1}^r\ii D_j=2.
$$
\end{lemma}

\bigskip
Let us return to the hypothesis of Theorem \ref{0} that
$X$ is a connected minimal set of $f$.
Given an invariant domain $U$, we have ${\rm Fix}(f)\cap U\neq\emptyset$
by Lemma \ref{9} and the Brouwer fixed point theorem applied
to the Carath\'eodory compactification $\hat U$.

\begin{lemma} \label{11}
The invariant domains are finite in number.
\end{lemma}

{\bf Proof}. Assume the invariant domains are infinite
and denote them by $U_i$ ($i=1,2,\cdots$).
Choose a fixed point $x_i$ from $U_i$. Then passing to a subsequence, $x_i$
converges to a point $x$ in $S^2$, which must be a fixed point
of $f$. 
If $x$ is contained in $X$, then $X$ has a fixed point, which contradicts the assumption.
Otherwise, $U_i$ coincide for large $i$.
A contradiction.
\qed

\bigskip
Choose a closed disc $D$
in $U$ which contains ${\rm Fix}(f)\cap U$
in its interior. Then $D$ is adapted and its index
$\ii D$ is independent of the choice of $D$.
Choose one of them and denote it by $D(U)$.

\begin{lemma} \label{12}
For any invariant domain $U$, the index $\ii D(U)$ is equal to 1.
\end{lemma}

{\bf Proof}. By Lemma \ref{9}, the Carath\'eodory
rotation number of $U$ is nonzero. On $\hat U$ the
region bounded by $\partial D(U)$ and $\PP(U)$ has no fixed point.
Thus one needs only compute the index of $\hat f$ w.\ r.\ t.\
the boundary curve $\PP(U)$.
\qed

\bigskip
Now let us conclude the proof of Theorem \ref{0}. Lemmata
\ref{10}, \ref{11} and \ref{12} clearly show that there are
exactly two invariant domains. 

For any $n>1$, the minimal set $X$ is minimal for $f^n$ since
it is connected. Applying the above result to $f^n$, one
can show that there is no more invariant domain of $f^n$.
Also the Carath\'eodory rotation number of an invariant domain
must be irrational, as is shown by applying Lemma \ref{9}
to the iterates of $f$. 

Finally that both Carath\'eodory rotation numbers
coincide follows from the main results of \cite{BG}.
The proof is complete.

\bigskip

Let us expose the Cartwright-Littlewood fixed point theorem.

\begin{theorem} \label{21}
Let $f$ be an orientation preserving homeomorphism of $S^2$.
Let $X$ be a continuum invariant by $f$. Assume $U=S^2\setminus X$
is connected. Then $f$ has a fixed point in $X$.
\end{theorem}

{\bf Proof}. 
Assume the contrary. 
If 
the Carath\'eodory rotation number of $U$ is
nonzero, then  Lemma \ref{12} shows that $\ii D(U)=1$.
If the rotation number is 0, 
Theorem \ref{5} says that the
homeomorphism $\hat f\vert_{\PP(U)}$ is  Morse-Smale, with $2n$
($n\geq1$) fixed points. Moreover the attractors
(resp.\ repellors) are attractors (resp.\ repellors) of the whole
map $\hat f$. In this case one can compute the index just following
the definition, with the result that $\ii D(U)=1-n$. 
Both cases contradicts
Lemma \ref{10}.
\qed

\section{Minimal continuum with wandering domain}

In \cite{Ha} a pathological $C^\infty$ diffeomorphism is constructed
which has a pseudo-circle $C$ as a minimal set.
See also \cite{He}. It is well known in continuum theory
that there are points $x$ in $C$
which are not accessible from both sides. Blowing up $x$, as well
as all the points of its orbit, we can construct a homeomorphism
which has a minimal continuum with wandering domain (see \cite{AO}).
Conversely if there are wandering domains whose domains $\{U_i\}$ satisfy that $\{\overline{U_i}\}$ is a null-sequence of mutually disjoint discs,
one can pinch each domain to a point, which characterize the complement of wandering domains (see \cite{BNW}).

\end{document}